arXiv:dg-ga/9512011v1 21 Dec 1995

# $L^2$-COHOMOLOGY OF GEOMETRICALLY INFINITE HYPERBOLIC 3-MANIFOLDS

JOHN LOTT

December 18, 1995

Abstract. We give results on the following questions about a topologically tame hyperbolic 3-manifold $M$ :
1. Does $M$ have nonzero $L^2$-harmonic 1-forms?
2. Does zero lie in the spectrum of the Laplacian acting on $\Lambda^1(M)/\mathrm{Ker}(d)$?

## 1. Introduction

Let $M$ be a complete oriented Riemannian manifold. A basic problem is to understand the spectrum of the Laplacian $\triangle_p$ acting on the square-integrable $p$-forms $\Lambda^p(M)$. In this paper we are concerned with the bottom of the spectrum. We address the following questions :
1. Does $M$ have nonzero $L^2$-harmonic $p$-forms?
2. Does zero lie in the spectrum of $\triangle_p$?

If $M$ is compact then Hodge theory tells us that questions 1 and 2 are equivalent and that the answer is "yes" if and only if $\mathrm{H}^p(M;\mathbb{C}) \neq 0$. In particular, the answer only depends on the topology of $M$.

If $M$ is noncompact then things are different. First, questions 1 and 2 are no longer equivalent - think of $M = \mathbb{R}$. Second, the answers to these questions no longer only depend on the topology of $M$. They depend on both the topology of $M$ and its asymptotic geometry in a subtle way which is not understood.

In this paper we look at the above questions for a class of Riemannian manifolds with interesting asymptotic geometry, namely connected hyperbolic 3-manifolds $M$ which are *topologically tame*, i.e. diffeomorphic to the interior of a compact 3-manifold with boundary. We review the relevant geometry of such manifolds in Section 3. Their ends can be characterized as cusps, flares and tubes. $M$ is called *geometrically finite* if its ends are all cusps or flares and *geometrically infinite* otherwise.

Using the Hodge decomposition, the square-integrable differential forms on $M$ can be split into $\mathrm{Ker}(\triangle_0)$, $\Lambda^0(M)/\mathrm{Ker}(d)$, $\mathrm{Ker}(\triangle_1)$ and $\Lambda^1(M)/\mathrm{Ker}(d)$. Hereafter we

Research supported by NSF grant DMS-9403652

assume that $M$ is noncompact. The only possible elements of $\mathrm{Ker}(\triangle_0)$ are constant functions and so if $\mathrm{vol}(M) < \infty$ then $\mathrm{Ker}(\triangle_0) = \mathbb{C}$, while if $\mathrm{vol}(M) = \infty$ then $\mathrm{Ker}(\triangle_0) = 0$. The next result of Canary tells what happens on $\Lambda^0(M)/\mathrm{Ker}(d)$ [4].

**Proposition 1.** *Zero lies in the spectrum of the Laplacian acting on* $\Lambda^0(M)/\mathrm{Ker}(d)$ *if and only if* $M$ *is geometrically infinite.*

Thus the spectrum of the Laplacian, acting on functions, is sensitive to whether $M$ has any tubular ends, but is not sensitive to the geometry of those ends. If $M$ is geometrically finite, Mazzeo and Phillips computed $\dim(\mathrm{Ker}(\triangle_1))$ and the essential spectrum of the Laplacian on $\Lambda^1(M)/\mathrm{Ker}(d)$ [13]. In particular, if $M$ is geometrically finite then zero always lies in the spectrum of the Laplacian acting on $\Lambda^1(M)/\mathrm{Ker}(d)$. One could ask whether there is a direct analogue of Canary's theorem for $\Lambda^1(M)/\mathrm{Ker}(d)$. However, the following example shows that such cannot be the case.

Let $S$ be a closed oriented surface of genus $g \geq 2$ and let $\phi \in \mathrm{Diff}(S)$ be an orientation-preserving pseudo-Anosov diffeomorphism of $S$. Thurston showed that the mapping torus $MT$ of $\phi$ has a hyperbolic metric [17, 22]. The corresponding cyclic cover $M$ of $MT$ is a geometrically infinite hyperbolic 3-manifold. In Section 4 we prove

**Proposition 2.** *Zero lies in the spectrum of the Laplacian acting on* $\Lambda^1(M)/\mathrm{Ker}(d)$ *if and only if* $\phi^* \in \mathrm{Aut}\left(\mathrm{H}^1(S;\mathbb{R})\right)$ *has an eigenvalue of norm one.*

It is known that any element of $\mathrm{Sp}(2g,\mathbb{Z})$ can occur as $\phi^*$ for some pseudo-Anosov diffeomorphism of $S$ [18]. Thus the result of Proposition 2 is not vacuous. It shows that the spectrum of the Laplacian, acting on 1-forms, is sensitive to the geometry of the tubular ends.

The manifolds considered in Proposition 2 are very special. The question arises how to extend Proposition 2 to general hyperbolic 3-manifolds $M$ of finite topological type. First, we dispose of the case when $M$ has zero injectivity radius. In Section 5 we prove

**Proposition 3.** *If* $\inf_{m \in M} \mathrm{inj}(m) = 0$ *then the essential spectrum of the Laplacian acting on* $\Lambda^1(M)/\mathrm{Ker}(d)$ *is* $[0, \infty)$.

We are left with the case of positive injectivity radius. There is an obvious problem in studying the spectrum of the Laplacian on $M$ in that we do not have an explicit description of the Riemannian metric of $M$. For example, even in the above case of a mapping torus, the hyperbolic metric on $MT$ is constructed by an iterative process. Our way of getting around this problem is to translate questions about the bottom of the spectrum into questions about the reduced and unreduced $L^2$-cohomology of $M$. It is much easier to compute the $L^2$-cohomologies of $M$ than to compute the spectral resolution of its Laplacian. Furthermore, the $L^2$-cohomologies of $M$ only depend on the biLipschitz diffeomorphism class of $M$. In our case we do know what $M$ looks like up to a biLipschitz diffeomorphism, thanks to the work of Minsky [15].

Let $M$ be a topologically tame hyperbolic 3-manifold with positive injectivity radius. We make the technical assumption that the ends of $M$ are incompressible. For brevity, we call such a hyperbolic 3-manifold *nice*. Minsky gave a length space which models the large-scale geometry of $M$. By a slight variation of his work, we construct a model manifold $\mathcal{M}$ which is biLipschitz diffeomorphic to $M$. The geometry of a tubular end $[0, \infty) \times S$ of $\mathcal{M}$ is given by a ray $\gamma$ in the Teichmüller space $\mathcal{T}_S$ of the surface $S$. The endpoint of $\gamma$, a point in Thurston's compactification of $\mathcal{T}_S$, is the *end invariant* of the tubular end. It is known that $M$ is determined up to isometry by its topology and its end invariants [15]. Hence the question is how exactly these determine the spectrum of the Laplacian.

Each point $\gamma(t)$ along the ray gives an inner product $\langle \cdot, \cdot \rangle_t$ on $\mathrm{H}^1(S; \mathbb{R})$. Let $\Gamma'(\mathrm{H}^1)$ be the Hilbert space of measurable maps $f : [0, \infty) \to \mathrm{H}^1(S; \mathbb{R})$ such that $\int_0^\infty \langle f(t), f(t) \rangle_t \, dt < \infty$. Put $\Gamma(\mathrm{H}^1) = \{ f \in \Gamma'(\mathrm{H}^1) : \partial_t f \in \Gamma'(\mathrm{H}^1) \}$. In Section 6 we prove

**Proposition 4.** *Let $M$ be a nice hyperbolic 3-manifold. Then zero is not in the spectrum of the Laplacian acting on $\Lambda^1(M)/\mathrm{Ker}(d)$ if and only if each end of $M$ is tubular and the corresponding operator $\partial_t : \Gamma(\mathrm{H}^1) \to \Gamma'(\mathrm{H}^1)$ has closed image.*

The next proposition gives a sufficient condition for $\partial_t$ to be onto. In Section 7 we prove

**Proposition 5.** *Suppose that there is a decomposition $\mathrm{H}^1(S; \mathbb{R}) = E_+ \oplus E_-$ and constants $a, c_+, c_- > 0$ such that for all $v_+ \in E_+$, $v_- \in E_-$ and $s_1 \geq s_2 \geq 0$,*

$$\|v_+\|_{s_1} \geq c_+ \, e^{a(s_1 - s_2)} \|v_+\|_{s_2}$$

*and*

$$\|v_-\|_{s_1} \leq c_- \, e^{-a(s_1 - s_2)} \|v_-\|_{s_2}.$$

*Then $\partial_t$ is onto.*

We also give a conjectural algorithm to determine directly from the end invariants whether or not zero lies in the spectrum of the Laplacian acting on $\Lambda^1(M)/\mathrm{Ker}(d)$, at least for most end invariants.

Finally, we give results on $\mathrm{Ker}(\triangle_1)$. In Section 6 we prove

**Proposition 6.** *If $M$ is a nice hyperbolic 3-manifold then $\dim(\mathrm{Ker}(\triangle_1)) < \infty$.*

Let $K$ be a compact submanifold of $M$ such that $M$ retracts onto $\mathrm{int}(K)$. Put $L_1 = \mathrm{Im}\left(\mathrm{H}^1(K; \mathbb{R}) \to \mathrm{H}^1(\partial K; \mathbb{R})\right)$. It is a Lagrangian subspace of $\mathrm{H}^1(\partial K; \mathbb{R})$. In Section 8 we prove

**Proposition 7.** *Let $M$ be a nice hyperbolic 3-manifold. Suppose that zero is not in the spectrum of the Laplacian acting on $\Lambda^1(M)/\mathrm{Ker}(d)$. For each end of $M$, consider the vector space* $\mathrm{Ker}\left(\partial_t : \Gamma(\mathrm{H}^1) \to \Gamma'(\mathrm{H}^1)\right)$. *Together, these give a Lagrangian subspace $L_2$ of* $\mathrm{H}^1(\partial K;\mathbb{R})$. *There is a short exact sequence*

$$0 \longrightarrow \mathrm{Im}\left(\mathrm{H}^1(K,\partial K;\mathbb{R}) \to \mathrm{H}^1(K;\mathbb{R})\right) \longrightarrow \mathrm{Ker}(\triangle_1) \longrightarrow L_1 \cap L_2 \to 0.$$

The organization of this paper is as follows. In Section 2 we define the reduced and unreduced $L^2$-cohomology groups and give their basic properties, along with their relation to the spectrum of the Laplacian. Some of these results are scattered throughout the literature, but we have tried to give a coherent presentation. In Section 3 we review the geometry of hyperbolic 3-manifolds and results of Minsky. In Section 4 we compute the reduced and unreduced $L^2$-cohomology groups of cyclic covers of general mapping tori. In Section 5 we consider hyperbolic 3-manifolds with vanishing injectivity radius. In Section 6 we describe the $L^2$-cohomology groups of tubular ends in terms of the operators $\partial_t : \Gamma(\mathrm{H}^1) \to \Gamma'(\mathrm{H}^1)$. In Section 7 we give sufficient conditions for the vanishing or nonvanishing of the unreduced $L^2$-cohomology groups of tubular ends. We also describe results of Zorich and their relation to spectral questions. In Section 8 we consider reduced $L^2$-cohomology groups of hyperbolic 3-manifolds.

I thank Josef Dodziuk and Rafe Mazzeo for discussions. I thank Yair Minsky and Anton Zorich for explanations of their work and for comments on parts of this paper. I especially thank Curt McMullen for many helpful conversations. I thank the IHES, the Max-Planck-Institut-Bonn and the Bonner Kaffeehaus for their hospitality.

## 2. $L^2$-cohomology

Let $M$ be an oriented Riemannian manifold which is geodesically complete except for a possible compact boundary. Consider the Hilbert space

$$\Lambda^p(M) = \{\text{square-integrable measurable}\, p-\text{forms on}\ M\} \tag{2.1}$$

and the subspace

$$\Omega^p(M) = \{\omega \in \Lambda^p(M) : d\omega \text{ is square-integrable on} \ \mathrm{int}(M)\}, \tag{2.2}$$

where $d\omega$ is initially interpreted in a distributional sense. There is a cochain complex

$$\ldots \xrightarrow{d_{p-1}} \Omega^p(M) \xrightarrow{d_p} \Omega^{p+1}(M) \xrightarrow{d_{p+1}} \ldots \tag{2.3}$$

One can check that $\mathrm{Ker}(d_p)$ is a closed subspace of $\Lambda^p(M)$.

**Definition 1.** *The $p$-th $L^2$-cohomology group of $M$ is* $\mathrm{H}^p_{(2)}(M) = \mathrm{Ker}(d_p)/\mathrm{Im}(d_{p-1})$. *The $p$-th reduced $L^2$-cohomology group of $M$ is* $\overline{\mathrm{H}}^p_{(2)}(M) = \mathrm{Ker}(d_p)/\overline{\mathrm{Im}(d_{p-1})}$, *a Hilbert space.*

We will sometimes call $\mathrm{H}^p_{(2)}(M)$ the $p$-th unreduced $L^2$-cohomology group. Let $M'$ be another manifold like $M$. Let $\Omega^*(M')$ be its cochain complex, with differential $d'$.

**Lemma 1.** *Suppose that there are linear maps*

$$\begin{aligned} i &: \Omega^*(M) \longrightarrow \Omega^*(M'), & K &: \Omega^*(M) \longrightarrow \Omega^{*-1}(M), \\ j &: \Omega^*(M') \longrightarrow \Omega^*(M), & K' &: \Omega^*(M') \longrightarrow \Omega^{*-1}(M') \end{aligned} \tag{2.4}$$

*such that*

$$\begin{aligned} i \circ d &= d' \circ i, & j \circ d' &= d \circ j, \\ I - j \circ i &= dK + Kd, & I - i \circ j &= d'K' + K'd'. \end{aligned} \tag{2.5}$$

*Then $j$ induces an isomorphism between $\mathrm{H}^*_{(2)}(M')$ and $\mathrm{H}^*_{(2)}(M)$. If $i$ and $j$ are continuous then $j$ also induces an isomorphism between $\overline{\mathrm{H}}^*_{(2)}(M')$ and $\overline{\mathrm{H}}^*_{(2)}(M)$.*

*Proof.* We leave the proof to the reader. □

The natural geometric invariance of $L^2$-cohomology turns out to be Lipschitz homotopy equivalence. We will only consider maps $f : M \to M'$ such that $f(\partial M) \subset \partial M'$.

**Definition 2.** *1. A map $f : M \to M'$ is said to be Lipschitz if $f$ is almost everywhere differentiable and there is a constant $C > 0$ such that for almost all $m \in M$ and all $v \in T_mM$, $|(df)_m v| \leq C|v|$.*
*2. Two Lipschitz maps $f_0 : M \to M'$ and $f_1 : M \to M'$ are Lipschitz-homotopic if there is a Lipschitz map $F : [0,1] \times M \to M'$ which restricts to $f_0$ and $f_1$ on the boundary.*
*3. Two Lipschitz maps $f : M \to M'$ and $g : M' \to M$ define a Lipschitz-homotopy equivalence between $M$ and $M'$ if $f \circ g$ and $g \circ f$ are Lipschitz-homotopic to the identity.*

A Lipschitz map $f : M \to M'$ induces maps $f^* : \mathrm{H}^*_{(2)}(M') \to \mathrm{H}^*_{(2)}(M)$ and $f^* : \overline{\mathrm{H}}^*_{(2)}(M') \to \overline{\mathrm{H}}^*_{(2)}(M)$.

**Proposition 8.** *If $f : M \to M'$ and $g : M' \to M$ define a Lipschitz-homotopy equivalence between $M$ and $M'$ then $f^*$ induces an isomorphism between $\mathrm{H}^*_{(2)}(M')$ and $\mathrm{H}^*_{(2)}(M)$, and between $\overline{\mathrm{H}}^*_{(2)}(M')$ and $\overline{\mathrm{H}}^*_{(2)}(M)$.*

*Proof.* The homotopy-equivalence gives continuous linear maps $i = g^*$, $j = f^*$, $K$ and $K'$ satisfying the hypotheses of Lemma 1. □

Let $\delta$ denote the formal $L^2$-adjoint of $d$. Let $*$ denote the Hodge duality operator. Let $b : \partial M \to M$ be the boundary inclusion. Let $\Lambda^*_\infty(M)$ denote the smooth compactly-supported forms on $M$. Note if $\omega \in \Lambda^*_\infty(M)$ then $b^*(\omega)$ may be nonzero. Define a sequence of inner products $\langle \cdot, \cdot \rangle_s$ on $\Lambda^*_\infty(M)$ for $s \in \mathbb{N}$ inductively by $\langle \cdot, \cdot \rangle_0 = \langle \cdot, \cdot \rangle_{L^2}$ and

$$\langle \omega_1, \omega_2 \rangle_{s+1} = \langle \omega_1, \omega_2 \rangle_s + \langle d\omega_1, d\omega_2 \rangle_s + \langle \delta\omega_1, \delta\omega_2 \rangle_s . \tag{2.6}$$

Define the Sobolev space $\mathcal{H}^*_s(M)$ to be the completion of $\Lambda^*_\infty(M)$ under $\langle \cdot, \cdot \rangle_s$.

The Laplace operator is $\triangle = \delta d + d\delta$. It is a self-adjoint operator on $\Lambda^*(M)$ with domain

$$\mathrm{Dom}(\triangle) = \{\omega \in \mathcal{H}^*_2(M) : b^*(*\omega) = b^*(*d\omega) = 0\} \tag{2.7}$$

and if $\dim(M) > 0$ then it is unbounded. If $\rho \in L^\infty([0,\infty))$ then $\rho(\triangle)$ is a bounded operator on $\Lambda^*(M)$. Let $\triangle_p$ be the restriction of $\triangle$ to $\Lambda^p(M)$. We have

$$\begin{aligned} \overline{\mathrm{H}}^p_{(2)}(M) &\cong \mathrm{Ker}(d_p) \cap (\mathrm{Im}(d_{p-1}))^\perp \\ &= \{\omega \in \Omega^p(M) : d\omega = \delta\omega = b^*(*\omega) = 0\} = \mathrm{Ker}(\triangle_p). \end{aligned} \tag{2.8}$$

By elliptic theory, $\mathrm{Ker}(\triangle_p)$ consists of smooth forms and so $\overline{\mathrm{H}}^*_{(2)}(M)$ can be computed using only smooth forms. We now show that the same is true for $\mathrm{H}^*_{(2)}(M)$. Put

$$\Omega^{p,\infty}(M) = \{\omega \in \Omega^p(M) : \omega \text{ is smooth}\}. \tag{2.9}$$

There is a complex

$$\ldots \xrightarrow{d_{p-1}} \Omega^{p,\infty}(M) \xrightarrow{d_p} \Omega^{p+1,\infty}(M) \xrightarrow{d_{p+1}} \ldots \tag{2.10}$$

**Proposition 9.** *The cohomology of the complex (2.10) is isomorphic to $\mathrm{H}^*_{(2)}(M)$.*

*Proof.* There is an obvious cochain map $i : \Omega^{*,\infty}(M) \to \Omega^*(M)$. Let $\eta \in C^\infty([0,\infty))$ be identically 1 on $[0,1]$ and identically 0 on $[2,\infty)$. Then $\eta(\triangle)$ is a smoothing operator and gives a cochain map $j : \Omega^*(M) \to \Omega^{*,\infty}(M)$. Define $\rho \in C^\infty([0,\infty))$ by $\rho(x) = \frac{1-\eta(x)}{x}$ and define $K : \Omega^*(M) \to \Omega^{*-1}(M)$ by $K = \delta\rho(\triangle)$. Then $I - ij = dK + Kd$ and similarly for $I - ji$. The proposition follows. □

We now show that the $L^2$-cohomology groups can be computed by means of standard elliptic complexes for manifolds with boundary.

For $s \in \mathbb{Z}$, there is a Hilbert cochain complex $\mathcal{D}_s(M)$ given by

$$0 \to \mathcal{H}^0_{s+dim(M)}(M) \to \mathcal{H}^1_{s+dim(M)-1}(M) \to \ldots \to \mathcal{H}^{dim(M)-1}_{s+1}(M) \to \mathcal{H}^{dim(M)}_s(M) \to 0, \tag{2.11}$$

where we implicitly truncate the complex when the Sobolev index becomes negative.

For fixed $p$, consider the Hilbert cochain complex $\mathcal{D}_{abs}(M)$, concentrated in degrees $p-1$, $p$ and $p+1$, given by

$$\begin{aligned} \mathcal{D}^{p-1}_{abs}(M) &= \{\omega \in \mathcal{H}^{p-1}_2(M) : b^*(*d\omega) = b^*(*\omega) = 0\}, \\ \mathcal{D}^p_{abs}(M) &= \{\omega \in \mathcal{H}^p_1(M) : b^*(*\omega) = 0\}, \\ \mathcal{D}^{p+1}_{abs}(M) &= \mathcal{H}^{p+1}_0(M). \end{aligned} \tag{2.12}$$

**Proposition 10.** *If $s \geq p + 1 - \dim(M)$ then the part of $\mathcal{D}_s(M)$ from degrees $p-1$ to $p+1$ is homotopy equivalent to $\mathcal{D}_{abs}(M)$.*

*Proof.* Let $\epsilon > 0$ be small enough that there is a coordinate function $t \in [0, 2\epsilon]$ near $\partial M$ such that $\partial_t$ is a unit length vector field whose flow generates unit speed geodesics which are normal to $\partial M$, and $\partial M$ corresponds to $t = 0$. Using these coordinates, a tubular neighborhood of $\partial M$ is diffeomorphic to $[0, 2\epsilon] \times \partial M$. Let $Y$ denote a copy of $M$ but with the product metric on $[0, 2\epsilon] \times \partial M$. The identity map gives a homotopy equivalence between $\mathcal{D}_s(M)$ and $\mathcal{D}_s(Y)$. Let $DY$ denote the double of $Y$ and let $\mathcal{D}_s^{even}(DY)$ be the complex of forms on $DY$ which are invariant under the $\mathbb{Z}_2$-involution on $DY$. There is an obvious inclusion $f : \mathcal{D}_s^{even}(DY) \to \mathcal{D}_s(Y)$. We now show that $\mathcal{D}_s(Y)$ and $\mathcal{D}_s^{even}(DY)$ are homotopy equivalent.

A differential form $\omega$ on $Y$ can be decomposed near the boundary as

$$\omega = \omega_1(t) + dt \wedge \omega_2(t), \tag{2.13}$$

where $\omega_1(t)$ and $\omega_2(t)$ are forms on $\partial M$. Let $\rho : [0, 2\epsilon] \to \mathbb{R}$ be a smooth bump function which is identically one near $t = 0$ and identically zero for $t \geq \epsilon$. Let $\widehat{\triangle}$ denote the Laplacian on $\partial M$. For $u > 0$, define the operator

$$R(u) = I - e^{-\frac{I+\widehat{\triangle}}{u^2}} \tag{2.14}$$

by the spectral theorem. For $\omega$ a form on $Y$, restrict $\omega$ to $[0, 2\epsilon] \times \partial M$ and put

$$(K\omega)(t) = \rho(t) \int_0^t R(u)\, \omega_2(u) du. \tag{2.15}$$

Then one can check that $K$ acts as a degree $-1$ map on both $\mathcal{D}_s(Y)$ and $\mathcal{D}_s^{even}(DY)$. If $\omega$ is a form on $Y$ then near $\partial M$,

$$\begin{aligned} \omega - (dK + Kd)\omega &= \omega_1(0) + ((I - R(t))\, \omega_1(t) + dt \wedge ((I - R(t))\, \omega_2(t) \\ &\quad + \int_0^t R'(u)\, \omega_1(u) du. \end{aligned} \tag{2.16}$$

One can check that $\omega - (dK + Kd)\omega$ extends by reflection to an element of $\mathcal{D}_s^{even}(DY)$. Thus we obtain a homotopy equivalence $f : \mathcal{D}_s^{even}(DY) \to \mathcal{D}_s(Y)$ and $g : \mathcal{D}_s(Y) \to \mathcal{D}_s^{even}(DY)$, where $f$ is the inclusion map and $g = I - (dK + Kd)$.

Next, as $s$ varies the complexes $\mathcal{D}_s^{even}(DY)$ are all isomorphic to each other by powers of $I + \triangle_{DY}$, at least in their common terms of definition. Thus we may consider the case $s = p + 1 - \dim(M)$. In this case, the part of $\mathcal{D}_s^{even}(DY)$ from $p - 1$ to $p + 1$ is the same as $\mathcal{D}_{abs}(Y)$.

Finally, we show that $\mathcal{D}_{abs}(M)$ is the same as $\mathcal{D}_{abs}(Y)$. Let us decompose a form $\omega$ on $M$ as in (2.13). Then the boundary condition for $\omega$ to belong to $\mathcal{D}^p_{abs}(M)$ is $\omega_2(0) = 0$ and the additional boundary condition for $\omega$ to belong to $\mathcal{D}^{p-1}_{abs}(M)$ is $\partial_t \omega_1(0) = 0$. These conditions determine the same spaces of forms whether one is on $M$ or $Y$. □

**Proposition 11.** *The reduced and unreduced p-th $L^2$-cohomology groups of $M$ are isomorphic to the reduced and unreduced p-th cohomology groups of the complex $\mathcal{D}_{abs}(M)$.*

*Proof.* For the reduced $L^2$-cohomology, the claim follows from (2.8). As the operator $(I+\triangle)^{-1/2}$ is an isomorphism from $\Lambda^p(M)$ to $\mathcal{D}^p_{abs}(M)$, it follows from Definition 1 that

$$\mathrm{H}^p_{(2)}(M) \cong \frac{\mathrm{Ker}(d) \text{ on } \mathcal{D}^p_{abs}(M)}{\mathrm{Im}(d) \text{ on } \{\omega \in \mathcal{H}_1^{p-1}(M) : b^*(*\omega) = 0, d\omega \in \mathcal{D}^p_{abs}(M)\}}. \tag{2.17}$$

The Hodge decomposition on $M$ is

(2.18)
$$\Lambda^*(M) = \mathrm{Ker}(\triangle_*) \oplus \overline{\mathrm{Im}(d) \text{ on } \mathcal{H}_1^{*-1}(M)} \oplus \overline{\mathrm{Im}(\delta) \text{ on } \{\omega \in \mathcal{H}_1^{*+1}(M) : b^*(*\omega) = 0\}}.$$

Projecting $\omega$ from (2.17) onto the last factor in (2.18), we may as well assume that $\delta\omega = 0$, showing that $\omega \in \mathcal{D}^{p-1}_{abs}(M)$. □

Let $i_p$ be the obvious surjection from $\mathrm{H}^p_{(2)}(M)$ to $\overline{\mathrm{H}}^p_{(2)}(M)$. We have $\mathrm{Ker}(i_{p+1}) = \overline{\mathrm{Im}(d_p)}/\mathrm{Im}(d_p)$. Thus $i_{p+1}$ is an isomorphism if and only if $\mathrm{Im}(d_p)$ is closed.

*For the rest of this section, we assume that $\partial M = \emptyset$.*

Let $K$ be a compact submanifold of $M$ with smooth boundary $\partial K$. Put $N = \overline{M - K}$.

**Proposition 12.** *We have that*
*1. The reduced $L^2$-cohomology at $p$ of $\mathcal{D}_{abs}(M)$ is finite-dimensional if and only if the reduced $L^2$-cohomology at $p$ of $\mathcal{D}_{abs}(N)$ is finite-dimensional.*
*2. The reduced $L^2$-cohomology at $p$ of $\mathcal{D}_{abs}(M)$ equals the unreduced $L^2$-cohomology if and only if the reduced $L^2$-cohomology at $p$ of $\mathcal{D}_{abs}(N)$ equals the unreduced $L^2$-cohomology.*

*Proof.* Let $Z$ be a small collaring of $\partial K$ in $M$, diffeomorphic to $[-1,1] \times \partial K$. Put $K' = K \cup Z$ and $N' = N \cup Z$. Then $K'$ is diffeomorphic to $K$ and $N'$ is diffeomorphic to $N$, with $K' \cap N' = Z$. Let $i_1 : K' \to M$, $i_2 : N' \to M$, $i_3 : Z \to K'$ and $i_4 : Z \to N'$ be the obvious embeddings. There is a short exact sequence

$$0 \longrightarrow \mathcal{D}_s(M) \stackrel{i_1^* \oplus i_2^*}{\longrightarrow} \mathcal{D}_s(K') \oplus \mathcal{D}_s(N') \stackrel{i_3^* - i_4^*}{\longrightarrow} \mathcal{D}_s(Z) \longrightarrow 0. \tag{2.19}$$

(Warning : One may be tempted to use the cohomology sequence of (2.19) to compare the reduced $L^2$-cohomology groups of $M$ and $N'$. However, this cohomology sequence need not be weakly exact if one does not make Fredholmness assumptions. We do not want to make such assumptions.)

Let us take smooth compatible triangulations of $K'$ and $Z$. Let $C^*(K')$ and $C^*(Z)$ denote the corresponding (finite-dimensional!) complexes of simplicial cochains. For $s$ large enough, it is known that $\mathcal{D}_s(K')$ is homotopy equivalent to $C^*(K')$ and $\mathcal{D}_s(Z)$ is homotopy equivalent to $C^*(Z)$ [6, 23]. Explicitly, the maps involved are integration $\int : \mathcal{D}_s(K') \to C^*(K')$ and the Whitney map $W : C^*(K') \to \mathcal{D}_s(K')$, and similarly for $Z$. We have commuting diagrams

$$\int : \quad \begin{array}{ccc} \mathcal{D}_s(K') \oplus \mathcal{D}_s(N') & \stackrel{q}{\rightarrow} & \mathcal{D}_s(Z) \\ \downarrow & & \downarrow \\ C^*(K') \oplus \mathcal{D}_s(N') & \stackrel{q'}{\rightarrow} & C^*(Z) \end{array} \tag{2.20}$$

and

$$W : \quad \begin{array}{ccc} \mathcal{D}_s(K') \oplus \mathcal{D}_s(N') & \stackrel{q}{\rightarrow} & \mathcal{D}_s(Z) \\ \uparrow & & \uparrow \\ C^*(K') \oplus \mathcal{D}_s(N') & \stackrel{q'}{\rightarrow} & C^*(Z), \end{array} \tag{2.21}$$

where $q = i_3^* - i_4^*$ and $q'(c, \omega) = i_3^* c - \int i_4^* \omega$. Also, the relevant homotopy operators form commuting diagrams. It follows that the complex $\mathrm{Ker}(q)$ is homotopy equivalent to $\mathrm{Ker}(q')$.

From (2.19), $\mathrm{Ker}(q) \cong \mathcal{D}_s(M)$. Note that $\mathrm{Ker}(q')$ has the Hilbert structure arising from its inclusion in $C^*(K') \oplus \mathcal{D}_s(N')$. Let $d$ denote the differential in $C^*(K') \oplus \mathcal{D}_s(N')$ and let $\delta$ denote its adjoint. Put

$$V = \{\omega \in C^*(K') \oplus \mathcal{D}_s(N') : d\omega = 0 \text{ and } \delta\omega \in (\mathrm{Ker}(q'))^\perp\}. \tag{2.22}$$

Using "harmonic representatives", we can identify the reduced $L^2$-cohomology of $\mathrm{Ker}(q')$ with $V \cap \mathrm{Ker}(q')$ and that of $C^*(K') \oplus \mathcal{D}_s(N')$ with $V \cap \mathrm{Ker}(\delta)$. As these are both of finite codimension in $V$ and $H^*(K')$ is finite-dimensional, it follows that the reduced $L^2$-cohomology of $\mathrm{Ker}(q')$ is finite-dimensional if and only if that of $\mathcal{D}_s(N')$ is finite-dimensional. As $N'$ is diffeomorphic to $N$ by a diffeomorphism which is an isometry outside of a compact region, the reduced $L^2$-cohomology of $\mathcal{D}_s(N')$ is isomorphic to that of $\mathcal{D}_s(N)$. Part 1 of the proposition now follows from Proposition 10.

Finally, let $c$ denote the differential in $\mathrm{Ker}(q')$. It follows from [9, Lemma 13.6.2] that $c$ has closed image if and only if $d$ has closed image. Part 2 of the proposition follows. □

**Proposition 13.**

$$0 \notin \sigma\,(\delta d \text{ on } \Lambda^p(M)/\mathrm{Ker}(d)) \iff i_{p+1} \text{ is an isomorphism.}$$

*Proof.* Suppose first that $\delta d$ has a bounded inverse on $\Lambda^p(M)/\mathrm{Ker}(d)$. Given $\mu \in \Lambda^p(M)$, let $\overline{\mu}$ denote its class in $\Lambda^p(M)/\mathrm{Ker}(d)$. Define an operator $S$ on smooth compactly-supported $(p+1)$-forms by $S(\omega) = d(\delta d)^{-1}\overline{\delta\omega}$. Then $S$ extends to a bounded operator on $\Lambda^{p+1}(M)$. Let $\{\eta_n\}_{n\in\mathbb{N}}$ be a sequence in $\Omega^p(M)$ such that $\lim_{n\to\infty} d\eta_n = \omega$ for some $\omega \in \Lambda^{p+1}(M)$. Then for each $n \in \mathbb{N}$, we have $d\eta_n = S(d\eta_n)$ and so $\omega = S(\omega)$. Thus $\omega \in \mathrm{Im}(d)$ and so $\mathrm{Im}(d)$ is closed.

Now suppose that $\delta d$ does not have a bounded inverse on $\Lambda^p(M)/\mathrm{Ker}(d)$. Then there is a sequence of positive numbers $r_1 > s_1 > r_2 > s_2 > \ldots$ tending towards zero and an orthonormal sequence $\{\eta_n\}_{n\in\mathbb{N}}$ in $\Lambda^p(M)/\mathrm{Ker}(d)$ such that with respect to the spectral projection $P$ of $\delta d$, $\eta_n \in \mathrm{Im}(P([r_n, s_n]))$. Put $\lambda_n = \|d\eta_n\|$. Then $\lim_{n\to\infty} \lambda_n = 0$. Let $\{c_n\}_{n\in\mathbb{N}}$ be a sequence in $\mathbb{R}^+$ such that $\sum_{n=1}^{\infty} c_n^2 = \infty$ and $\sum_{n=1}^{\infty} c_n\lambda_n < \infty$. Put $\omega = \sum_{n=1}^{\infty} c_n d\eta_n$. Then $\omega \in \overline{\mathrm{Im}(d)}$. Suppose that $\omega = d\mu$ for some $\mu \in \Omega^p(M)$. By the spectral theorem, we must have $\overline{\mu} = \sum_{n=1}^{\infty} c_n\eta_n$. However, this is not square-integrable. Thus $\mathrm{Im}(d)$ is not closed. $\square$

We recall the notion of the *essential spectrum* of an operator. Let $T$ be a densely-defined self-adjoint operator on a Hilbert space $H$. Then $\sigma_{ess}(T)$ is a closed subset of the spectrum $\sigma(T)$ with the property that $\lambda \in \sigma_{ess}(T) \iff 0 \in \sigma_{ess}(T - \lambda I)$. Let $P$ be the spectral projection of $T$. Then $\sigma_{ess}(T)$ has the following equivalent characterizations [10].

**Proposition 14.** $0 \in \sigma_{ess}(T)$ *if and only if any of the following conditions hold :*
*1.* $\dim(\mathrm{Ker}(T)) = \infty$ *or* $\mathrm{Im}(T)$ *is not closed.*
*2. There is a bounded sequence* $\{u_n\}_{n\in\mathbb{N}}$ *in* $\mathrm{Dom}(T)$ *such that* $\lim_{n\to\infty} \|Tu_n\| = 0$*, but* $\{u_n\}_{n\in\mathbb{N}}$ *does not have a convergent subsequence.*
*3. There is an orthonormal sequence* $\{u_n\}_{n\in\mathbb{N}}$ *in* $\mathrm{Dom}(T)$ *such that* $\lim_{n\to\infty} \|Tu_n\| = 0$.
*4. For all* $\epsilon > 0$, $\dim(\mathrm{Im}(P([-\epsilon, \epsilon]))) = \infty$.
*5.* $\dim(\mathrm{Ker}(T)) = \infty$ *or* $0$ *is not isolated in* $\sigma(T)$.

In particular, if $\mathrm{Ker}(T) = 0$ then $0 \in \sigma_{ess}(T) \iff 0 \in \sigma(T)$.

**Corollary 1.** *Let* $M$ *and* $M'$ *be complete oriented Riemannian manifolds. Suppose that there are compact submanifolds* $K \subset M$ *and* $K' \subset M'$ *such that* $M - K$ *is isometric to* $M' - K'$*. Then*
*1.* $0 \in \sigma_{ess}(\triangle_p \text{ on } \mathrm{Ker}(\triangle_p(M))) \iff 0 \in \sigma_{ess}(\triangle_p \text{ on } \mathrm{Ker}(\triangle_p(M')))$.
*2.* $0 \in \sigma_{ess}(\triangle_p \text{ on } \Lambda^p(M)/\mathrm{Ker}(d)) \iff 0 \in \sigma_{ess}(\triangle_p \text{ on } \Lambda^p(M')/\mathrm{Ker}(d))$.
*3.* $0 \in \sigma_{ess}(\triangle_p \text{ on } \Lambda^p(M)) \iff 0 \in \sigma_{ess}(\triangle_p \text{ on } \Lambda^p(M'))$.

*Proof.* 1. As $\triangle_p$ acts on $\mathrm{Ker}(\triangle_p(M))$ as the zero operator, Proposition 14.1 says that 0 lies in $\sigma_{ess}(\triangle_p \text{ on } \mathrm{Ker}(\triangle_p(M)))$ if and only if $\dim(\mathrm{Ker}(\triangle_p(M))) = \infty$. The claim follows from (2.8) and Proposition 12.1.
2. As $\triangle_p$ acts on $\Lambda^p(M)/\mathrm{Ker}(d)$ as $\delta d$, the claim follows from Propositions 12.2 and 13.
3. This is now a consequence of the Hodge decomposition. $\square$

*Remark* : Corollary 1.3 is well-known. It is a consequence of [8, Prop. 2.1], the proof of which is for functions but extends to differential forms. We will need the more refined statements of Corollary 1.1, 1.2, which take into account the Hodge decomposition of forms on $M$ and $M'$. The proof of [8, Prop. 2.1], which involves multiplication by a compactly supported function, does not extend to this case. Consequently, we have given an independent proof. I thank Jozef Dodziuk for correspondence on these questions.

## 3. Hyperbolic 3-Manifolds

For background on hyperbolic 3-manifolds, we refer to [1, 20, 21]. Let $M = H^3/\Gamma$ be a complete connected oriented hyperbolic 3-manifold with finitely generated fundamental group $\Gamma$. We assume that $\Gamma$ is nonabelian, as the abelian case can be easily handled separately. The sphere at infinity of $H^3$ breaks up as the union $S^2 = \Lambda \cup \Omega$ of the *limit set* and the *domain of discontinuity*, on which $\Gamma$ acts freely. The *convex core* of $M$ is $C(M) = \mathrm{hull}(\Lambda)/\Gamma$. The quotient $\Omega/\Gamma$ is a finite union of connected Riemann surfaces, each of which is diffeomorphic to the complement of a finite number of points in a closed connected Riemann surface. Put $\overline{M} = (H^3 \cup \Omega)/\Gamma$.

There is a constant $\mu$, the Margulis constant, such that if $\epsilon < \mu$ and

$$M_{thin}(\epsilon) = \{m \in M : \mathrm{inj}(m) < \epsilon\}$$

then each connected component of $M_{thin}(\epsilon)$ is either
1. A rank-two cusp, diffeomorphic to $(0,\infty) \times T^2$,
2. A rank-one cusp, diffeomorphic to $(0,\infty) \times (-1,1) \times S^1$, or
3. A tubular neighborhood of a short geodesic loop in $M$, diffeomorphic to $S^1 \times D^2$.

Let $M^0(\epsilon)$ be $M$ with the cusps in $M_{thin}(\epsilon)$ removed. There is a notion of an *end* $E$ of $M^0(\epsilon)$ and of $E$ being contained in an open set $U \subset M^0(\epsilon)$ [2]. An end $E$ is said to be *geometrically finite* if it is contained in an open set $U$ such that $U \cap C(M) = \emptyset$. If $E$ is geometrically finite then it is associated to a connected component of $\Omega/\Gamma$. The complex structure on that component is called the *end invariant* of $E$. $M$ is said to be *geometrically finite* if all of the ends of $M^0(\epsilon)$ are geometrically finite and *geometrically infinite* otherwise.

If $M$ is geometrically finite then there is a pair $(X, P)$, where $X$ is a compact 3-manifold and $P$ is a compact submanifold of $\partial X$, with the property that $M$ is diffeomorphic to $\mathrm{int}(X)$ and $\overline{M}$ is diffeomorphic to $X - P$. The *parabolic locus* $P$ is a disjoint union $T \cup A$ of surfaces, where $T$ is a disjoint union of 2-tori, one for each rank-two cusp of $M$, and $A$ is a disjoint union of annuli, one for each rank-one cusp of $M$.

The reduced $L^2$-cohomology and essential spectrum of geometrically finite hyperbolic manifolds were studied in [13]. When specialized to three dimensions, the results are as follows. If $M$ has finite volume then $\overline{\mathrm{H}}^0_{(2)}(M) \cong \mathbb{C}$ and if $M$ has infinite volume then $\overline{\mathrm{H}}^0_{(2)}(M) = 0$. The first reduced $L^2$-cohomology group of $M$ is given by

$$\overline{\mathrm{H}}^1_{(2)}(M) \cong \mathrm{Im}\left(\mathrm{H}^1(X, \partial X - \mathrm{int}(A)) \longrightarrow \mathrm{H}^1(X, \partial X - (T \cup \mathrm{int}(A))\right). \tag{3.1}$$

The essential spectrum of $\triangle$ is

(3.2)

| | $\Lambda^0/\mathrm{Ker}(d)$ | $\Lambda^1/\mathrm{Ker}(d)$ |
|---|---|---|
| $M$ compact | $\emptyset$ | $\emptyset$ |
| $M$ noncompact | $[1, \infty)$ | $[0, \infty)$ |

We no longer assume that $M$ is geometrically finite. The fact that $\Gamma$ is finitely-generated implies that $M$ is homotopy-equivalent to a compact 3-manifold [19]. It is an open conjecture, which has been proved in many cases, that $M$ must be *topologically tame*, i.e. diffeomorphic to the interior of a compact 3-manifold. We assume hereafter that $M$ is topologically tame. There is again a pair $(X, P)$, where $X$ is a compact 3-manifold and $P$ is a compact submanifold of $\partial X$, with the properties that
1. $M$ is diffeomorphic to $\mathrm{int}(X)$.
2. $P$ is a union of tori and annuli, one for each cuspidal component of $M_{thin}(\epsilon)$.
3. The ends of $M^0(\epsilon)$ are in one-to-one correspondence with the connected components of $\partial X - P$.

An end $E$ of $M^0(\epsilon)$ is called *simply degenerate* if it is contained in an open set $U \subset M^0(\epsilon)$ homeomorphic to $(0, \infty) \times S$ for some compact surface $S$, and there is a sequence of finite-area hyperbolic surfaces in $U$, each homotopic to $\{0\} \times \mathrm{int}(S)$, such that the sequence exits the end; see [5] for the precise definition. It is known that $M$ is *geometrically tame*, meaning that every end of $M^0(\epsilon)$ is either geometrically finite or simply degenerate [2, 5]. A simply degenerate end $E$ comes equipped with a certain geodesic lamination on the surface $\mathrm{int}(S)$, known as its *end invariant*. Let $\mathcal{E}$ denote the collection of all end invariants of $M$. Thurston conjectured that $M$ is determined up to isometry by the topology of $(X, P)$, along with $\mathcal{E}$ [21]. We remark that if the triple $(X, P, \mathcal{E})$ satisfies certain topological conditions then it does arise from some hyperbolic 3-manifold [16].

Canary showed that if $M$ is geometrically infinite then $C(M)$, an infinite volume submanifold of $M$, can be exhausted by compact submanifolds whose boundary areas are uniformly bounded above [5]. As he pointed out in [4], it then follows from Buser's theorem that zero lies in the spectrum of the Laplacian acting on $L^2$-functions on $M$.

Suppose that there is a constant $c > 0$ such that for all $m \in M$, $\mathrm{inj}(m) > c$. Then $P = \emptyset$ and any simply degenerate end of $M$ is contained in an open set of the form $(0, \infty) \times S$ for some closed oriented surface $S$. Suppose in addition that the ends of $M$ are incompressible, or equivalently, that $\Gamma$ does not decompose as a nontrivial free product. In this case, Minsky showed that Thurston's isometry conjecture is true [15]. To do so, Minsky first constructed a model Riemannian manifold $\mathcal{M}$, based on the topology of $M$ and its end invariants, which approximates the geometry of $M$.

More precisely, he showed that there is a map $f : \mathcal{M} \to M$ which is homotopic to a homeomorphism, with the property that the lift $\widetilde{f} : \widetilde{\mathcal{M}} \to H^3$ is a coarse quasi-isometry. The Riemannian metric on $\mathcal{M}$ is constructed as follows. It is enough to first specify the Riemannian metric on the ends of $\mathcal{M}$ and then extend it arbitrarily to the rest of $\mathcal{M}$. If $U = (0, \infty) \times S$ contains a geometrically finite end of $M$, let $d\rho^2$ be the hyperbolic metric on the corresponding connected component of $\Omega/\Gamma$. Then the model metric on the associated end of $\mathcal{M}$ is $dt^2 + \cosh^2(t) d\rho^2$.

To describe the model metric for a simply degenerate end of $M$, we first need some notation. For a closed oriented surface $S$ of genus $g \geq 2$, let $\mathcal{H}_S$ denote the space of hyperbolic metrics on $S$, let $\mathrm{Diff}_S$ denote the orientation-preserving diffeomorphisms of $S$ and let $\mathrm{Diff}_{0,S}$ denote those isotopic to the identity. The Teichmüller space $\mathcal{T}_S$ can be identified with $\mathcal{H}_S/\mathrm{Diff}_{0,S}$ and the moduli space $\mathrm{Mod}_S$ can be identified with $\mathcal{H}_S/\mathrm{Diff}_S$. Note that $\mathrm{Mod}_S$ is an orbifold. There is a quotient map $\pi : \mathcal{T}_S \to \mathrm{Mod}_S$. The universal Teichmüller curve $\widehat{\mathcal{T}}_S$ is $\mathcal{H}_S \times_{Diff_{0,S}} S$. It is the total space of a fiber bundle $p_T : \widehat{\mathcal{T}}_S \to \mathcal{T}_S$ with fiber $S$ and inherits an obvious family of hyperbolic metrics on its fibers. The universal curve over the moduli space is $\widehat{\mathrm{Mod}}_S = \mathcal{H}_S \times_{Diff_S} S$. It is the total space of an orbifold fiber bundle $p_M : \widehat{\mathrm{Mod}}_S \to \mathrm{Mod}_S$ with fiber $S$ and again inherits a family of hyperbolic metrics on its fibers. Let us choose an arbitrary smooth horizontal distribution on the fiber bundle $\widehat{\mathrm{Mod}}_S$, meaning a collection of subspaces $T^H\widehat{\mathrm{Mod}}_S \subset T\widehat{\mathrm{Mod}}_S$ such that $T\widehat{\mathrm{Mod}}_S = T^H\widehat{\mathrm{Mod}}_S \oplus \mathrm{Ker}(dp_M)$. (Everything here is interpreted in an orbifold sense.) There is a lifted horizontal distribution $T^H\widehat{\mathcal{T}}_S$ on $\widehat{\mathcal{T}}_S$.

If $U = (0, \infty) \times S$ contains a simply degenerate end of $M$, fix an initial hyperbolic metric $d\rho^2(0)$ on $\{0\} \times S$ and let $S_0$ be the corresponding Riemann surface. Let $\mathrm{H}^0(S_0; K^2)$ denote the space of holomorphic quadratic differentials on $S_0$. It is a complex vector space of dimension $-\frac{3}{2}\chi(S)$. The ending lamination $L$ is equivalent to the vertical foliation of some $\Phi \in \mathrm{H}^0(S_0; K^2)$. Then $\Phi$ generates a Teichmüller ray $\gamma : [0, \infty) \to \mathcal{T}_S$ starting from $\gamma(0) = [S_0]$. The endpoint of $\gamma$ corresponds to $L$, viewed as a point in Thurston's compactification of $\mathcal{T}_S$. As the injectivity radius of $M$ is bounded below by a positive number, [15, Theorem 5.5] implies that the projected ray $\pi \circ \gamma$ lies in a compact region of $\mathrm{Mod}_S$.

Using the hyperbolic metrics on the fibers of $\widehat{\mathcal{T}}_S$, the horizontal distribution $T^H\widehat{\mathcal{T}}_S$ and the metric $dt^2$ on $[0, \infty)$, there is an induced Riemannian metric on $p_T^{-1}(\gamma)$. In terms of the trivialization $p_T^{-1}(\gamma) \cong [0, \infty) \times S$ coming from $T^H\widehat{\mathcal{T}}_S$, we can write this metric as $dt^2 + d\rho^2(t)$, where for each $t \in [0, \infty)$, $d\rho^2(t) \in \mathcal{H}_S$ projects to $\gamma(t) \in \mathcal{T}_S$. This is the model metric on the associated end of $\mathcal{M}$. Because of the precompactness of $\pi \circ \gamma$, the biLipschitz class of the model metric is independent of the choice of $T^H\widehat{\mathrm{Mod}}_S$.

We will need to know that $\mathcal{M}$ approximates $M$ in a slightly better way than that given in [15]. Curt McMullen explained to me how the next statement follows from the results of [15].

**Proposition 15.** *There is a biLipschitz homeomorphism between $\mathcal{M}$ and $M$.*

*Proof.* It is enough to just construct biLipschitz homeomorphisms between open sets containing the ends of $\mathcal{M}$ and $M$. For a geometrically finite end, this follows from [15, Theorem 5.2]. Let $E$ be a simply degenerate end of $M$ contained in a neighborhood $U = (0, \infty) \times S$. Let $\mathcal{U} = (0, \infty) \times S$ contain the corresponding end in $\mathcal{M}$. Let $\gamma$ be the Teichmüller ray described above. Minsky constructed a sequence $\{g_n : S \rightarrow U\}_{n \in \mathbb{N}}$ of pleated surfaces in $U$ with the properties [15, Theorem 5.5] that
1. for each $n \in \mathbb{N}$, $g_n(S)$ is homotopic in $U$ to $\{0\} \times S$
2. the sequence $\{g_n(S)\}_{n \in \mathbb{N}}$ exits the end
3. there is a constant $T > 0$ such that for each $n \in \mathbb{N}$, the Teichmüller class of the induced hyperbolic metric $\rho_n \in \mathcal{H}_S$ lies within a Teichmüller distance $T$ from $\gamma(n)$.

After precomposing the $g_n$'s with appropriate elements of $\mathrm{Diff}_{0,S}$, we may assume that neighboring $\rho_n$'s are uniformly close in the sense that there is a $K > 0$ such that for all $n \in \mathbb{N}$, the identity map $\mathrm{Id} : (S, \rho_n) \rightarrow (S, \rho_{n+1})$ is a $K$-biLipschitz homeomorphism.

For each $n \in \mathbb{N}$, we can find an embedded surface $h_n : S \rightarrow U$ such that $h_n(S)$ lies within some distance $D$ from $g_n(S)$ and the induced hyperbolic metric $\rho_n'$ is $K'$-biLipschitz to $\rho_n$ for some $K' > 0$. As the injectivity radius of $M$ is bounded below by a positive constant, we can use compactness in the geometric topology [1, Chapter E], [14, Section 4] to argue that the surfaces can be chosen so that $D$ and $K'$ are uniform with respect to $n$. Next, we can find constants $0 < a < A$ and a uniformly spaced subsequence $\{n_k\}_{k \in \mathbb{N}}$ of $\mathbb{N}$ so that the consecutive surfaces $\{h_{n_k}(S), h_{n_k+1}(S)\}$ are spaced at least distance $a$ apart and no more than distance $A$. Using property 2. above, we can assume that the surfaces $\{h_{n_k}(S)\}$ are topologically consecutive in the sense that $h_{n_k}(S)$ separates $h_{n_{k-1}}(S)$ from $h_{n_{k+1}}(S)$. Let $U_k$ be the part of $U$ enclosed by $h_{n_k}(S)$ and $h_{n_k+1}(S)$. Let $\mathcal{U}_k$ be the submanifold $[n_k, n_{k+1}] \times S$ in $\mathcal{U}$. For each $k \in \mathbb{N}$, there is a diffeomorphism $\phi_k : \mathcal{U}_k \rightarrow U_k$ which restricts to $\{h_{n_k}, h_{n_k+1}\}$ on $\partial \mathcal{U}_k$. Again using compactness in the geometric topology, we can choose the diffeomorphisms $\{\phi_k\}_{k \in \mathbb{N}}$ so that there is a constant $K'' > 0$ such that for all $k \in \mathbb{N}$, $\phi_k$ is a $K''$-biLipschitz homeomorphism. The desired biLipschitz homeomorphism $f : \mathcal{U} \rightarrow U$ is given by stacking together the $\phi_k$'s. □

*Remark* : Minsky used singular Euclidean metrics on $S$ instead of hyperbolic metrics, but the difference is minor. We use the horizontal distribution on $\widehat{\mathrm{Mod}}_S$ to give a lifting of $\gamma$ to $\mathcal{H}_S$ such that the lifts of nearby points on $\gamma$ are uniformly close in $\mathcal{H}_S$. This is similar to [15, p. 562-563].

## 4. Mapping Tori

Let $F$ be a smooth closed oriented manifold. Let $\phi \in \mathrm{Diff}(F)$ be an orientation-preserving diffeomorphism of $F$. The mapping torus of $\phi$ is the manifold

$$MT = ([0,1] \times F) \,/ \sim \tag{4.1}$$

where the equivalence relation is $(0,s) \sim (1,\phi(s))$.

Projection on the first factor gives a fibering $\pi : MT \to S^1$. Let $M$ be the associated cyclic cover of $MT$. Let $\phi_p^* \in \mathrm{Aut}(\mathrm{H}^p(F,\mathbb{R}))$ be the map on cohomology coming from $\phi$.

**Proposition 16.** *1.* $\overline{\mathrm{H}}^*_{(2)}(M) = 0$.
*2.* $0 \in \sigma\,(\delta d \,\mathrm{on}\, \Lambda^p(M)/\mathrm{Ker}(d)) \iff \phi_p^*$ *has an eigenvalue of norm one.*

*Proof.* Let $\gamma$ denote a generator of the group of covering transformations on $M$, the one which maps to $t \to t+1$ on $\mathbb{R}$. For $\lambda \in U(1)$, put

$$\Lambda^p_\lambda(M) = \{\mathrm{measurable}\, p - \mathrm{forms}\, \omega \,\mathrm{on}\, M : \gamma^*\omega = \lambda\omega\}. \tag{4.2}$$

Let $V_\lambda$ be the flat complex line bundle on $S^1$ with holonomy $\lambda$ and put $E_\lambda = \pi^* V_\lambda$. Then

$$\Lambda^p_\lambda(M) \cong \Lambda^p(MT; E_\lambda), \tag{4.3}$$

the $p$-forms on $MT$ with value in the flat vector bundle $E_\lambda$. It follows from Fourier analysis that there is a direct integral decomposition

$$\Lambda^p(M) = \int_{U(1)} \Lambda^p(MT; E_\lambda)\, d\lambda. \tag{4.4}$$

Furthermore, the decomposition (4.4) commutes with the Laplacians. It follows from Floquet theory that $\overline{\mathrm{H}}^p_{(2)}(M) \neq 0$ if and only if $H^p(MT; E_\lambda) \neq 0$ for all $\lambda \in U(1)$ and $0 \in \sigma(\triangle_p(M))$ if and only if $H^p(MT; E_\lambda) \neq 0$ for some $\lambda \in U(1)$; see [11] for details.

There is a Wang exact sequence

(4.5)
$$\ldots \to H^{p-1}(F) \overset{I-\lambda^{-1}\phi^*_{p-1}}{\to} H^{p-1}(F) \to H^p(MT; E_\lambda) \to H^p(F) \overset{I-\lambda^{-1}\phi^*_p}{\to} H^p(F) \to \ldots$$

This gives the short exact sequence

$$0 \longrightarrow \mathrm{Coker}(I - \lambda^{-1}\phi^*_{p-1}) \longrightarrow H^p(MT; E_\lambda) \longrightarrow \mathrm{Ker}(I - \lambda^{-1}\phi^*_p) \longrightarrow 0. \tag{4.6}$$

As there is only a finite number of $\lambda \in U(1)$ such that $\mathrm{Coker}(I - \lambda^{-1}\phi^*_{p-1}) \neq 0$ or $\mathrm{Ker}(I - \lambda^{-1}\phi^*_p) \neq 0$, part 1 of the proposition follows.

The Hodge decomposition of $\Lambda^p(M)$ now gives

$$\Lambda^p(M) = \overline{\mathrm{Im}(d \,\mathrm{on}\, \Lambda^{p-1}(M)/\mathrm{Ker}(d))} \oplus \Lambda^p(M)/\mathrm{Ker}(d). \tag{4.7}$$

Correspondingly, we have

(4.8)
$$0 \in \sigma\left(\triangle_p \text{ on } \overline{\mathrm{Im}(d \text{ on } \Lambda^{p-1}(M)/\mathrm{Ker}(d))}\right) \iff \mathrm{Coker}(I - \lambda^{-1}\phi^*_{p-1}) \neq 0,$$
$$0 \in \sigma\left(\triangle_p \text{ on } \Lambda^p(M)/\mathrm{Ker}(d)\right) \iff \mathrm{Ker}(I - \lambda^{-1}\phi^*_p) \neq 0$$

for some $\lambda \in U(1)$. The proposition follows. □

*Remark* : A different proof of Proposition 16 follows from Appendix A of the preprint version of [12]. This material was left out in the printed version.

Now let $F$ be a closed oriented surface $S$ of genus $g \geq 2$. Let $\phi$ be an orientation-preserving pseudo-Anosov diffeomorphism of $S$. Thurston showed that the mapping torus $MT$ has a hyperbolic structure [17, 22]. Furthermore, the hyperbolic structure on $MT$ is unique up to isometry. The cyclic cover $M$ has the pullback hyperbolic structure.

**Corollary 2.** $0 \in \sigma\left(\delta d \text{ on } \Lambda^1(M)/\mathrm{Ker}(d)\right) \iff \phi^*_1$ *has an eigenvalue of norm one.*

## 5. Zero Injectivity Radius

**Proposition 17.** *Let $M$ be a complete hyperbolic 3-manifold with* $\inf_{m \in M} \mathrm{inj}(m) = 0$. *Then* $\sigma_{ess}\left(\delta d \text{ on } \Lambda^1(M)/\mathrm{Ker}(d)\right) = [0, \infty)$.

*Proof.* If $M_{thin}(\epsilon)$ contains cusps then the proposition follows from the characterization of the essential spectrum of cusps in [13]. Otherwise, $M_{thin}(\epsilon)$ must contain a sequence of tubular neighborhoods $\{T_n\}_{n \in \mathbb{N}}$ of short geodesic loops $\{\gamma_n\}_{n \in \mathbb{N}}$ whose lengths $\{l(\gamma_n)\}_{n \in \mathbb{N}}$ tend towards zero. It is known that the radius $R_n$ of $T_n$ goes like $R_n \sim \frac{1}{2} \log\left(\frac{1}{l(\gamma_n)}\right)$ [7]. As $n$ increases, the geometry of $T_n$ approaches that of a rank-two cusp and so the claim of the proposition is not surprising.

Fix $n$ for a moment. We use Fermi coordinates $(r, t, \theta)$ for $T_n$ as in [7], where $0 \leq r \leq R_n$ is the distance to $\gamma_n$, $t$ is the arc-length along $\gamma_n$ and $\theta$ is the angular coordinate in the normal disk bundle to $\gamma_n$. Consider a 1-form $\omega$ on $T_n$ given in coordinates by $\omega = g(r)dt$, where $g \in C^\infty_0(0, R_n)$. One can check [7] that $\delta\omega = 0$,

$$\langle \omega, \omega \rangle = 2\pi l(\gamma_n) \int_0^{R_n} |g(r)|^2 \tanh(r) dr \tag{5.1}$$

and

$$\delta d\omega = -\frac{1}{\tanh(r)} \left(\tanh(r)\, g'\right)'. \tag{5.2}$$

Furthermore, $\omega \in \mathrm{Im}(\delta)$ if

$$\int_0^{R_n} g(r) \tanh(r) dr = 0, \tag{5.3}$$

or equivalently, if $\langle \omega, dt \rangle = 0$.

Let $\phi \in C_0^\infty((0,1))$ be a positive function satisfying $\int_0^1 \phi^2(r)dr = 1$. For $k \in \mathbb{R}$, define

$$g_{n,k}(r) = \frac{1}{\sqrt{2\pi l(\gamma_n) R_n}} e^{ikr} \phi(r/R_n) \tag{5.4}$$

and $\omega_{n,k} = g_{n,k}(r)dt$. We now fix $k \neq 0$. Put

$$\begin{aligned} c_{n,k} = \frac{\langle \omega_{n,k}, dt \rangle}{\langle \omega_{n,0}, dt \rangle} &= \frac{\int_0^{R_n} g_{n,k}(r) \ \tanh(r)dr}{\int_0^{R_n} g_{n,0} \ \tanh(r)dr} \\ &= \frac{\int_0^{R_n} e^{ikr}\phi(r/R_n) \ \tanh(r)dr}{\int_0^{R_n} \phi(r/R_n) \ \tanh(r)dr} \\ &= \frac{\int_0^1 e^{ikR_n s}\phi(s) \ \tanh(R_n s)ds}{\int_0^1 \phi(s) \ \tanh(R_n s)ds}. \end{aligned} \tag{5.5}$$

By the Riemann-Lebesgue theorem, $\lim_{n\to\infty} c_{n,k} = 0$. Put $\omega'_{n,k} = \omega_{n,k} - c_{n,k}\omega_{n,0}$. By construction, $\omega'_{n,k} \in \mathrm{Im}(\delta)$. We have

$$\begin{aligned} \|\omega'_{n,k}\|^2 &= 2\pi l(\gamma_n) \int_0^{R_n} |g_{n,k}(r) - c_{n,k} g_{n,0}(r)|^2 \tanh(r)dr \\ &= \frac{1}{R_n} \int_0^{R_n} \left| e^{ikr} - c_{n,k} \right|^2 \phi^2(r/R_n) \tanh(r)dr \\ &= \int_0^1 \left| e^{ikR_n s} - c_{n,k} \right|^2 \phi^2(s) \tanh(R_n s)ds \\ &= \int_0^1 \left(1 + c_{n,k}^2 - 2c_{n,k} \cos(kR_n s)\right) \phi^2(s) \tanh(R_n s)ds. \end{aligned} \tag{5.6}$$

Thus $\lim_{n\to\infty} \|\omega'_{n,k}\| = 1$.

Similarly, one can check that $\lim_{n\to\infty} \|(\delta d - k^2)\omega_{n,k}\| = 0$. It follows that

$$\lim_{n\to\infty} \|(\delta d - k^2)\omega'_{n,k}\| = \lim_{n\to\infty} \|(\delta d - k^2)\omega_{n,k} + k^2 c_{n.k}\omega_{n,0}\| = 0. \tag{5.7}$$

Since the $\omega'_{n,k}$'s are supported in the disjoint tubes $\{T_n\}_{n\in\mathbb{N}}$, they are mutually orthogonal. It follows from Proposition 14.3 that $k^2 \in \sigma_{ess}\,(\delta d \text{ on } \mathrm{Ker}(\delta) \subset \Lambda^1(M))$. As $\sigma_{ess}$ is a closed subset of $\mathbb{R}$, the proposition follows. □

*Remark* : There are hyperbolic 3-manifolds diffeomorphic to $\mathbb{R} \times S$, where $S$ is a closed oriented surface of genus $g \geq 2$, having zero injectivity radius [3].

## 6. Reduction to an ODE

Let $M$ be a topologically tame complete connected oriented hyperbolic 3-manifold. In this section, we are interested in whether zero lies in the spectrum of $\delta d$ acting on $\Lambda^1(M)/\mathrm{Ker}(d)$. If $M$ has zero injectivity radius then by Section 5, the essential spectrum of $\delta d$ acting on $\Lambda^1(M)/\mathrm{Ker}(d)$ is $[0, \infty)$. Therefore, we assume that $M$ has positive injectivity radius $c$.

We can take the constant $\epsilon$ in Section 3 less than $c$, so that $M_{thin}(\epsilon) = \emptyset$ and $M^0(\epsilon) = M$. By Section 2, it is enough to study the spectrum of the Laplacian on the ends of $M$. If $M$ has a geometrically finite end then it follows from [13] that $\sigma_{ess}\,(\delta d \text{ on } \Lambda^1(M)/\mathrm{Ker}(d)) = [0, \infty)$. Therefore, we assume that $M$ does not have any geometrically finite ends. By Section 3, all of the ends of $M$ are simply degenerate.

In order to apply Minsky's results, we make the further assumption that the ends of $M$ are incompressible. Recall from Section 3 that $\mathcal{M}$ is a certain Riemannian manifold which models $M$. By Propositions 8 and 15, $\mathrm{H}^2_{(2)}(M) \cong \mathrm{H}^2_{(2)}(\mathcal{M})$. Consider a single end of $\mathcal{M}$ which is contained in an open set $\mathcal{U} = (0, \infty) \times S$, where $S$ is a closed oriented surface. Our strategy will be to compute $\mathrm{H}^2_{(2)}(\overline{\mathcal{U}})$. Recall that $\mathcal{U}$ has the metric $dt^2 + d\rho^2(t)$, where $d\rho^2(t)$ is a hyperbolic metric on $S$ which projects to $\gamma(t) \in \mathcal{T}_S$.

For each $t \in [0, \infty)$, $\partial_t(d\rho^2(t))$ is a covariant 2-tensor on $S$. For any $k \in \mathbb{N}$, let $\|\partial_t(d\rho^2(t))\|_k$ denote its Sobolev $k$-norm with respect to $d\rho^2(t)$.

**Proposition 18.** *$\|\partial_t(d\rho^2(t))\|_k$ is uniformly bounded in $t$.*

*Proof.* As $d\rho^2(t)$ is a hyperbolic metric for all $t \in [0, \infty)$, it follows that

$$\partial_t(d\rho^2(t)) = \mathcal{L}_{V(t)} d\rho^2(t) + H(t) \tag{6.1}$$

where $V(t)$ is a vector field on $S$, $\mathcal{L}$ is the Lie derivative and $H(t)$ is a covariant 2-tensor on $S$ satisfying

$$\sum_\mu H^\mu{}_\mu(t) = 0, \quad \sum_\mu \nabla^\mu H_{\mu\nu}(t) = 0. \tag{6.2}$$

Equivalently, $H(t) = \mathrm{Re}(Q(t))$ where $Q(t) \in \mathrm{H}^0(S; K^2)$, $S$ having the complex structure induced from $d\rho^2(t)$. Let $z$ be a local holomorphic coordinate on $S$, write $d\rho^2(t)$ locally as $g_{z\overline{z}}\, dz\overline{dz}$ and write $Q(t)$ locally as $Q_{zz} dz^2$. The Beltrami differential corresponding to $Q(t)$ is

$$\mu = g^{z\overline{z}} Q_{zz}\, \frac{dz}{d\overline{z}}. \tag{6.3}$$

Given $\Phi = \Phi_{zz} dz^2 \in \mathrm{H}^0(S; K^2)$, put

$$\|\Phi\|_1 = \frac{i}{2} \int_S |\Phi_{zz}| \, dz \wedge d\overline{z}, \tag{6.4}$$
$$\|\Phi\|_2 = \frac{i}{2} \int_S g^{z\overline{z}} \Phi_{zz} \overline{\Phi_{zz}} \, dz \wedge d\overline{z}.$$

As $\gamma$ is a Teichmüller ray, the infinitesimal Teichmüller norm of $\gamma'(t)$ is

$$1 = \sup_{\{\Phi : \|\Phi\|_1 = 1\}} \left| \mathrm{Re} \left( \frac{i}{2} \int_S \Phi_{zz} \mu^z{}_{\overline{z}} \, dz \wedge d\overline{z} \right) \right| \tag{6.5}$$
$$= \sup_{\{\Phi : \|\Phi\|_1 = 1\}} \left| \langle \Phi, Q \rangle_2 \right| .$$

We now use that fact that $\pi \circ \gamma$ is precompact in $\mathrm{Mod}_S$. From the construction of $d\rho^2(t)$ using the horizontal distribution $T^H \widehat{\mathrm{Mod}}_S$, it follows that $\|\mathcal{L}_{V(t)} d\rho^2(t)\|_k$ is uniformly bounded in $t$. From (6.5), it follows that for fixed $t \in [0, \infty)$, $Q(t)$ lies in a compact subset of $\mathrm{H}^0(S; K^2)$ and hence one has a bound on $\|\mathrm{Re}(Q(t))\|_k$. Again using the precompactness of $\pi \circ \gamma$, it follows that $\|\mathrm{Re}(Q(t))\|_k$ is uniformly bounded in $t$. The proposition follows. □

For each $t \in [0, \infty)$, the vector space $\mathrm{H}^1(S; \mathbb{R})$ inherits a inner product $\langle \cdot, \cdot \rangle_t$ which can be described in two equivalent ways :
1. Given $h \in \mathrm{H}^1(S; \mathbb{R})$, let $\omega \in \Lambda^1(S)$ be its harmonic representative. Then

$$\langle h, h \rangle_t = \langle \omega, \omega \rangle_{d\rho^2(t)} = \int_S \omega \wedge *_t \omega. \tag{6.6}$$

2. Using the complex structure on $S$ coming from $\gamma(t)$, we can write $\mathrm{H}^1(S; \mathbb{R}) \otimes \mathbb{C} = \mathrm{H}^{1,0}(S) \oplus \mathrm{H}^{0,1}(S)$. Given $h \in \mathrm{H}^1(S; \mathbb{R})$, write $h = \frac{1}{2}(\rho + \overline{\rho})$ with $\rho \in \mathrm{H}^{1,0}(S)$. Then

$$\langle h, h \rangle_t = \frac{i}{2} \int_S \rho \wedge \overline{\rho}. \tag{6.7}$$

Let $\mathcal{H}^1(t)$ be the vector space of harmonic 1-forms on $S$, with respect to the metric $d\rho^2(t)$. Let $\Pi(t) : \Lambda^1(S; \mathbb{R}) \to \mathcal{H}^1(t)$ be the harmonic projection operator. Fix a set $\{C_i\}$ of closed $L^2$ 1-currents on $S$ whose homology representatives $\{[C_i]\}$ form a basis of $\mathrm{H}_1(S; \mathbb{R})$. Let $\{\tau^i\}$ be the dual basis of $\mathrm{H}^1(S; \mathbb{R})$. Define $\int_C : \Omega^1(S) \to \mathrm{H}^1(S; \mathbb{R})$ by

$$\int_C \omega = \sum_i \langle C_i, \omega \rangle \, \tau^i. \tag{6.8}$$

Then $\int_C$ restricts to an isomorphism $\int_C : \mathcal{H}^1(t) \to \mathrm{H}^1(S; \mathbb{R})$.

Let $\mathrm{H}^1$ be the vector bundle on $[0, \infty)$ whose fiber over $t \in [0, \infty)$ is isomorphic to $\mathrm{H}^1(S; \mathbb{R})$, with the inner product $\langle \cdot, \cdot \rangle_t$. Let $\mathcal{H}^1$ be the vector bundle on $[0, \infty)$ whose fiber over $t \in [0, \infty)$ is isomorphic to $\mathcal{H}^1(t)$.

**Definition 3.** *We define the following spaces.*
*1. Let* $\Gamma(\mathrm{H}^1)$ *be the vector space of* $L^2$*-sections* $\eta_1$ *of* $\mathrm{H}^1$ *such that* $\partial_t \eta_1 \in L^2([0,\infty);\mathrm{H}^1)$*.*
*2. Let* $\Gamma(\mathcal{H}^1)$ *be the vector space of* $L^2$*-sections* $\eta_1$ *of* $\mathcal{H}^1$ *such that* $(\Pi(t)\partial_t)\eta_1 \in L^2([0,\infty);\mathcal{H}^1)$*.*
*3. Let* $\Gamma'(\mathrm{H}^1)$ *be the vector space of* $L^2$*-sections of* $\mathrm{H}^1$*.*
*4. Let* $\Gamma'(\mathcal{H}^1)$ *be the vector space of* $L^2$*-sections of* $\mathcal{H}^1$*.*

There is an operator $\partial_t$ acting on $\Gamma(\mathrm{H}^1)$. Similarly, there is an operator $\Pi(t)\partial_t$ acting on $\Gamma(\mathcal{H}^1)$.

**Lemma 2.** *There is a commutative diagram*

$$\begin{array}{ccc} \Gamma(\mathcal{H}^1) & \xrightarrow{\Pi(t)\partial_t} & \Gamma'(\mathcal{H}^1) \\ \int_C \downarrow & & \int_C \downarrow \\ \Gamma(\mathrm{H}^1) & \xrightarrow{\partial_t} & \Gamma'(\mathrm{H}^1). \end{array} \tag{6.9}$$

*Proof.* Given $\omega \in \Gamma(\mathcal{H}^1)$,

$$\begin{aligned} \partial_t \int_C \omega &= \partial_t \sum_i \langle C_i, \omega\rangle \tau^i = \sum_i \langle C_i, \partial_t \omega\rangle \tau^i \\ &= \sum_i \langle C_i, \Pi(t)\partial_t\omega\rangle \tau^i = \int_C \Pi(t)\partial_t\omega. \end{aligned} \tag{6.10}$$

The lemma follows. □

Thus $\partial_t$, acting on $\Gamma(\mathrm{H}^1)$, is essentially the same as $\Pi(t)\partial_t$, acting on $\Gamma(\mathcal{H}^1)$.

Given $\eta \in \Omega^1(\overline{\mathcal{U}})$, write

$$\eta = \eta_1(t) + dt \wedge \eta_0(t), \tag{6.11}$$

where for almost all $t \in [0,\infty)$, $\eta_0(t) \in \Lambda^0(S)$ and $\eta_1(t) \in \Lambda^1(S)$. Define $\int_C : \Omega^1(\overline{\mathcal{U}}) \to \Gamma(\mathrm{H}^1)$ by

$$\left(\int_C \eta\right)(t) = \int_C \eta_1(t). \tag{6.12}$$

For each $t$, $\int_C : \Omega^1(S) \to \mathrm{H}^1(S;\mathbb{R})$ is bounded with respect to the inner products coming from $d\rho^2(t)$. By precompactness of $\pi \circ \gamma$, we can find a bound which is uniform in $t$. Thus $\int_C : \Omega^1(\overline{\mathcal{U}}) \to \Gamma(\mathrm{H}^1)$ is a bounded operator.

Given $\omega \in \Omega^2(\overline{\mathcal{U}})$, write

$$\omega = \omega_2(t) + dt \wedge \omega_1(t), \tag{6.13}$$

where for almost all $t \in [0,\infty)$, $\omega_1(t) \in \Lambda^1(S)$ and $\omega_2(t) \in \Lambda^2(S)$. Define $\int_C : \Omega^2(\overline{\mathcal{U}}) \to \Gamma'(\mathrm{H}^1)$ by

$$\left(\int_C \omega\right)(t) = \int_C \omega_1(t). \tag{6.14}$$

It is also a bounded operator.

**Lemma 3.** *There is a commutative diagram*

$$\begin{array}{ccccc} \Omega^1(\overline{\mathcal{U}}) & \xrightarrow{d_1} & \Omega^2(\overline{\mathcal{U}}) & \xrightarrow{d_2} & \Omega^3(\overline{\mathcal{U}}) \\ \int_C \downarrow & & \int_C \downarrow & & \downarrow \\ \Gamma(\mathrm{H}^1) & \xrightarrow{\partial_t} & \Gamma'(\mathrm{H}^1) & \longrightarrow & 0 \end{array} \tag{6.15}$$

*Proof.* Given $\eta$ as in (6.11),

$$\int_C d_1\eta = \int_C (\partial_t \eta_1 - d_S \eta_0) = \partial_t \int_C \eta_1 = \partial_t \int_C \eta. \tag{6.16}$$

The lemma follows. □

**Definition 4.** *Let* $J : \mathrm{H}^2_{(2)}(\overline{\mathcal{U}}) \to \Gamma'(\mathrm{H}^1)/\mathrm{Im}(\partial_t)$ *and* $\overline{J} : \overline{\mathrm{H}}^2_{(2)}(\overline{\mathcal{U}}) \to \Gamma'(\mathrm{H}^1)/\overline{\mathrm{Im}(\partial_t)}$ *be the maps induced from* (6.15)*.*

**Proposition 19.** *$J$ and $\overline{J}$ are isomorphisms.*

*Proof.* Let $\widehat{d}$ denote exterior differentiation on $\overline{\mathcal{U}}$ and let $d$ denote exterior differentiation on $S$. The condition for $\omega \in \Omega^2(\overline{\mathcal{U}})$ to be closed is

$$\partial_t \omega_2(t) = d\omega_1(t). \tag{6.17}$$

The equation $\omega = \widehat{d}\eta$ is equivalent to

$$\begin{aligned} \omega_1 &= \partial_t \eta_1 - d\eta_0, \\ \omega_2 &= d\eta_1. \end{aligned} \tag{6.18}$$

To see that $J$ and $\overline{J}$ are onto, take $h \in \Gamma'(\mathrm{H}^1)$. Choose harmonic representatives $\omega_1(t)$ for $h(t)$. Then $\omega = dt \wedge \omega_1(t)$ is closed and $\int_C \omega = h$.

We now show that $J$ is injective. By Proposition 9, we may assume that all forms considered are smooth. Suppose that $\omega \in \Omega^1(\mathcal{U})$ satisfies (6.17) and $\int_C \omega \in \mathrm{Im}(\partial_t)$. We want to find $\eta \in \Omega^1(\overline{\mathcal{U}})$ satisfying (6.18). We first show that we can eliminate $\omega_2$.

For any $t \in [0, \infty)$, let $[\omega_2(t)] \in \mathrm{H}^2(S; \mathbb{R}) \cong \mathbb{R}$ denote the de Rham cohomology class of $\omega_2(t)$. By (6.17), it is constant in $t$. By the Hodge decomposition, we have

$$\|\omega_2(t)\|^2 \geq [\omega_2(t)]^2 \mathrm{Area}_t(S) = -2\pi [\omega_2(t)]^2 \chi(S). \tag{6.19}$$

As $\omega \in \Lambda^2(\overline{\mathcal{U}})$, we must have $[\omega_2(t)] = 0$.

Let $G(t)$ be the Green's operator on $\Lambda^*(S)$, with respect to the metric $d\rho^2(t)$. It satisfies

$$\triangle(t)G(t) = G(t)\triangle(t) = I - \Pi(t), \quad \Pi(t)G(t) = G(t)\Pi(t) = 0. \tag{6.20}$$

We abbreviate $\delta(t)G(t)\omega_2(t) \in \Lambda^1(\overline{\mathcal{U}})$ by $\delta G \omega_2$.

**Lemma 4.** *$\widehat{d}(\delta G \omega_2)$ is square-integrable on $\overline{\mathcal{U}}$.*

*Proof.* The Hilbert space structure on $\Lambda^*(S)$ depends on $t$, but the underlying topological vector space is independent of $t$. Thus it makes sense to differentiate operators on $\Lambda^*(S)$ with respect to $t$. As $\pi \circ \gamma$ is precompact in $\mathrm{Mod}_S$, for any $s \in \mathbb{N}$, $G(t)$ is uniformly bounded in $t$ as an operator between the Sobolev spaces $\mathcal{H}^*_s(S)$ and $\mathcal{H}^*_{s+2}(S)$.

We have

$$\begin{aligned} \widehat{d}(\delta G \omega_2) &= d\delta G \omega_2 + dt \wedge \partial_t(\delta G \omega_2) \\ &= \omega_2 + dt \wedge ([\partial_t, \delta G]\omega_2 + \delta G \partial_t \omega_2) \\ &= \omega_2 + dt \wedge ([\partial_t, \delta G]\omega_2 + \delta G d\omega_1). \end{aligned} \tag{6.21}$$

As $\omega \in \Lambda^2(\overline{\mathcal{U}})$, it follows that $\omega_2$ and $\delta G d\omega_1$ are square-integrable on $\overline{\mathcal{U}}$. We must show that $[\partial_t, \delta G]\omega_2$ is square-integrable on $\overline{\mathcal{U}}$.

We have

$$[\partial_t, \delta G] = (\partial_t \delta) G + \delta(\partial_t G). \tag{6.22}$$

Acting on $\Lambda^2(S)$,

$$\delta = *^{-1}\, d\, * \tag{6.23}$$

and hence

$$\partial_t \delta = [\delta, *^{-1}(\partial_t *)]. \tag{6.24}$$

From Proposition 18, the Sobolev 1-norm of $*^{-1}(\partial_t *)$ is uniformly bounded in $t$. It follows that the operator $(\partial_t \delta)G$ is uniformly bounded in $t$.

Differentiating (6.20) with respect to $t$ gives

$$\partial_t G = -(\partial_t \Pi) G - G(\partial_t \Pi) - G(\partial_t \triangle) G. \tag{6.25}$$

Acting on $\phi \in \Lambda^2(S)$,

$$\Pi(t)\phi = \frac{\int_S \phi}{\mathrm{Area}_t(S)}\, d\mathrm{vol}_S(t) = -\frac{\int_S \phi}{2\pi\chi(S)}\, d\mathrm{vol}_S(t). \tag{6.26}$$

Thus

$$(\partial_t \Pi)\phi = -\frac{\int_S \phi}{2\pi\chi(S)}\, \partial_t d\mathrm{vol}_S. \tag{6.27}$$

Furthermore,

$$\partial_t \triangle = d(\partial_t \delta) + (\partial_t \delta) d. \tag{6.28}$$

Arguing as before, it follows that the operator $\delta(\partial_t G)$ is uniformly bounded in $t$. $\square$

Subtracting $\widehat{d}(\delta G\omega_2)$ from $\omega$, we may assume that $\omega_2 = 0$. We are left with a square-integrable $\omega_1$ satisfying

$$d\omega_1(t) = 0, \quad \int_C \omega_1 \in \mathrm{Im}(\partial_t). \tag{6.29}$$

We want to find square-integrable $\eta_0$ and $\eta_1$ satisfying

$$d\eta_1(t) = 0, \quad \partial_t\eta_1(t) - d\eta_0(t) = \omega_1(t). \tag{6.30}$$

Choose $h \in \Gamma(\mathrm{H}^1)$ such that $\int_C \omega_1 = \partial_t h$. Let $\eta_1(t)$ be the harmonic representative of $h(t)$. By definition,

$$d\eta_1(t) = \delta(t)\eta_1(t) = \triangle(t)\eta_1(t) = 0. \tag{6.31}$$

Then

$$d(\partial_t\eta_1 - \omega_1) = \partial_t d\eta_1 - d\omega_1 = 0. \tag{6.32}$$

Also

$$\int_C (\partial_t\eta_1 - \omega_1) = \partial_t h - \int_C \omega_1 = 0, \tag{6.33}$$

implying that

$$\Pi(t)(\partial_t\eta_1 - \omega_1) = 0. \tag{6.34}$$

Put

$$\eta_0 = \delta G(\partial_t\eta_1 - \omega_1). \tag{6.35}$$

Then

$$\begin{aligned} d\eta_0 &= d\delta G(\partial_t\eta_1 - \omega_1) = Gd\delta(\partial_t\eta_1 - \omega_1) \\ &= G(d\delta + \delta d)(\partial_t\eta_1 - \omega_1) = (I - \Pi(t))(\partial_t\eta_1 - \omega_1) \\ &= \partial_t\eta_1 - \omega_1. \end{aligned} \tag{6.36}$$

Thus equation (6.30) is satisfied.

It remains to show that $\eta$ is square-integrable. By construction, this is the case for $\eta_1$. We first show that $\partial_t\eta_1 - \omega_1$ is square-integrable. By assumption, this is the case for $\omega_1$. It follows from Lemma 2 that it is also the case for $(\Pi(t)\partial_t)\eta_1$. Differentiating (6.31) with respect to $t$ gives

$$(\partial_t\triangle)\eta_1 + \triangle(\partial_t\eta_1) = 0 \tag{6.37}$$

and so

$$\begin{aligned} (I - \Pi(t))(\partial_t\eta_1) &= G\triangle(\partial_t\eta_1) = -G(\partial_t\triangle)\eta_1 \\ &= -G(d(\partial_t\delta) + (\partial_t\delta)d)\eta_1 = -Gd(\partial_t\delta)\eta_1. \end{aligned} \tag{6.38}$$

Using the precompactness of $\pi \circ \gamma$ and arguing as in the proof of Lemma 4, it follows that $(I - \Pi(t))(\partial_t \eta_1)$ is square-integrable. Thus $\partial_t \eta_1$ is square-integrable.

We now have that $\partial_t \eta_1 - \omega_1$ is square-integrable. Again using the precompactness of $\pi \circ \gamma$ and arguing as in the proof of Lemma 4, it follows that from (6.35) that $\eta_0$ is square-integrable.

Finally, we show that $\overline{J}$ is injective. Suppose that $\omega \in \mathrm{Ker}(\widehat{d})$ and $\int_C \omega \in \overline{\mathrm{Im}(\partial_t)}$. We must show that $\omega \in \overline{\mathrm{Im}(\widehat{d})}$. As before, we may eliminate $\omega_2$. Thus we have a square-integrable $\omega_1$ and a sequence $\{h_n\}_{n \in \mathbb{N}}$ in $\Gamma(\mathrm{H}^1)$ satisfying

$$d\omega_1(t) = 0, \quad \int_C \omega_1 = \lim_n \partial_t h_n. \tag{6.39}$$

Let $\eta_{1,n}(t)$ be the harmonic representative of $h_n(t)$. Then

$$d(\partial_t \eta_{1,n} - \omega_1) = \partial_t d\eta_{1,n} - d\omega_1 = 0. \tag{6.40}$$

Put

$$\eta_{0,n} = \delta G(\partial_t \eta_{1,n} - \omega_1). \tag{6.41}$$

As before, $\eta_n = \eta_{1,n}(t) + dt \wedge \eta_{0,n}(t)$ is square-integrable. We have

$$\begin{aligned} d\eta_{0,n} &= d\delta G(\partial_t \eta_{1,n} - \omega_1) = G d\delta(\partial_t \eta_{1,n} - \omega_1) \\ &= G(d\delta + \delta d)(\partial_t \eta_{1,n} - \omega_1) = (I - \Pi(t))(\partial_t \eta_{1,n} - \omega_1). \end{aligned} \tag{6.42}$$

Thus if we can show that $\lim_n \Pi(t)(\partial_t \eta_{1,n} - \omega_1) = 0$, it will follow that $\omega = \lim_n \widehat{d}\eta_n$.

Now

$$\int_C \Pi(t)(\partial_t \eta_{1,n} - \omega_1) = \int_C (\partial_t \eta_{1,n} - \omega_1) = \partial_t h_n - \int_C \omega_1. \tag{6.43}$$

By precompactness of $\pi \circ \gamma$, the map $\int_C : \Gamma'(\mathcal{H}^1) \to \Gamma'(\mathrm{H}^1)$ is a Hilbert space isomorphism. Using (6.39), the proposition follows. □

**Proposition 20.** $\overline{\mathrm{H}}^2_{(2)}(\overline{\mathcal{U}}) = 0$.

*Proof.* From Proposition 19, $\overline{\mathrm{H}}^2_{(2)}(\overline{\mathcal{U}}) \cong \mathrm{Im}(\partial_t)^\perp \subset \Gamma'(\mathrm{H}^1)$. Using the inner product on $\Gamma'(\mathrm{H}^1)$, we can identify it with its dual space $\Gamma'(\mathrm{H}_1)$. Given $k \in \mathrm{Im}(\partial_t)^\perp$, let $\widetilde{k}$ be the corresponding element of $\Gamma'(\mathrm{H}_1)$. Let $h \in \Gamma'(\mathrm{H}^1)$ be smooth with compact support in $(0, \infty)$. As

$$0 = \langle k, \partial_t h \rangle = \int_0^\infty \left( \widetilde{k}(t), \partial_t h(t) \right) dt \tag{6.44}$$

holds for all such $h$, $\widetilde{k}(t)$ must be constant in $t$. Letting $h$ now have compact support in $[0, \infty)$, (6.44) gives that $\widetilde{k} = 0$. Hence $k = 0$. □

Using an $L^2$ 0-current and a closed $L^2$ 2-current on $S$, we can extend (6.15) to a commutative diagram

(6.45)
$$\begin{array}{ccccccc} \Omega^0(\overline{\mathcal{U}}) & \stackrel{d_0}{\to} & \Omega^1(\overline{\mathcal{U}}) & \stackrel{d_1}{\to} & \Omega^2(\overline{\mathcal{U}}) & \stackrel{d_2}{\to} & \Omega^3(\overline{\mathcal{U}}) \\ \downarrow & & \downarrow & & \downarrow & & \downarrow \\ \Omega^0([0,\infty)) & \stackrel{\partial_t}{\to} & \Omega^1([0,\infty)) \oplus \Gamma(\mathrm{H}^1) & \stackrel{\partial_t}{\to} & \Gamma'(\mathrm{H}^1) \oplus \Omega^0([0,\infty)) & \stackrel{\partial_t}{\to} & \Omega^1([0,\infty)). \end{array}$$

**Proposition 21.** *We have*

$$\overline{\mathrm{H}}^1_{(2)}(\overline{\mathcal{U}}) \cong \mathrm{Ker}\left(\partial_t : \Gamma(\mathrm{H}^1) \to \Gamma'(\mathrm{H}^1)\right) \tag{6.46}$$

*and*

(6.47)
$$\mathrm{H}^1_{(2)}(\overline{\mathcal{U}}) \cong \mathrm{Ker}\left(\partial_t : \Gamma(\mathrm{H}^1) \to \Gamma'(\mathrm{H}^1)\right) \oplus \frac{\Omega^1([0,\infty))}{\mathrm{Im}\left(\partial_t : \Omega^0([0,\infty)) \to \Omega^1([0,\infty))\right)}.$$

*In particular,* $\overline{\mathrm{H}}^1_{(2)}(\overline{\mathcal{U}})$ *is isomorphic to a subspace of* $\mathrm{H}^1(S;\mathbb{R})$.

*Proof.* The proof is similar to that of Proposition 19. We omit the details. □

**Corollary 3.** *Let $N$ be a connected oriented Riemannian 3-manifold. Suppose that there is a compact submanifold $K$ of $N$ such that each connected component $C_i$ of $N-K$ is isometric to a geometrically finite or simply degenerate end $E_i$ of a hyperbolic 3-manifold $M_i$. Suppose that $N$ has injectivity radius bounded below by a positive constant and each $M_i$ has incompressible ends. Then*
*1.* $\dim(\mathrm{Ker}(\triangle_1(N))) < \infty$
*2.* $0 \notin \sigma\left(\delta d \text{ on } \Lambda^1(N)/\mathrm{Ker}(d)\right)$ *if and only if each end of $N$ is geometrically infinite and the corresponding operator* $\partial_t : \Gamma(\mathrm{H}^1) \to \Gamma'(\mathrm{H}^1)$ *has closed image.*

*Proof.* Equation (2.8) and Propositions 8, 12, 13 and 15 imply that it is enough to verify the claims for the corresponding ends of the model manifolds $\mathcal{M}_i$.
1. If an end is geometrically finite, the claim follows from (3.1). If an end is geometrically infinite, the claim follows from Proposition 21.
2. If an end is geometrically finite, the claim follows from (3.2). If an end is geometrically infinite, the claim follows from Proposition 19. □

*Remark* : Corollary 3.1 is not an immediate consequence of the fact that $N$ has finite topological type. For example, the analogous statement for hyperbolic surfaces would be false.

## 7. Unreduced $L^2$-Cohomology

In Section 6 we reduced the problem of computing the $L^2$-cohomologies of an end of $M$ to that of computing the image of the operator $\partial_t$ on $\Gamma(\mathrm{H}^1)$. The inner product $\langle \cdot, \cdot \rangle_t$ on $\Gamma(\mathrm{H}^1)$ is determined by the Teichmüller geodesic $\gamma$. The question now arises as to how $\langle \cdot, \cdot \rangle_t$ depends on $t$.

**Example 1 :** Consider the mapping torus $MT$ discussed at the end of Section 4, whose fiber is a closed oriented surface $S$ of genus $g \geq 2$ and whose monodromy is an orientation-preserving pseudo-Anosov diffeomorphism $\phi$ of $S$. Let $\{d\rho^2(t)\}_{t \in \mathbb{R}}$ be a smooth curve in $\mathcal{H}_S$ such that for all $t \in \mathbb{R}$, $d\rho^2(t) = \phi^*(d\rho^2(t+1))$. Such a curve can be constructed by choosing an arbitrary $d\rho^2(0) \in \mathcal{H}_S$, choosing an arbitrary path $\{d\rho^2(t)\}_{t \in [0,1]}$ from $d\rho^2(0)$ to $(\phi^{-1})^*(d\rho^2(0))$ and then perturbing the path near the ends if necessary so that it extends to give $\{d\rho^2(t)\}_{t \in \mathbb{R}}$. The metric $dt^2 + d\rho^2(t)$ on $\mathbb{R} \times S$ descends to a metric on $MT$. Thus $dt^2 + d\rho^2(t)$ serves as a model metric for the hyperbolic metric on the cyclic cover $M$.

As $\phi^*$ acts symplectically on $\mathrm{H}^1(S;\mathbb{R})$, there is a decomposition

$$\mathrm{H}^1(S;\mathbb{R}) = E_0 \oplus \bigoplus_{i=1}^{k} (E_i \oplus E_{-i}) \tag{7.1}$$

and positive numbers

$$\lambda_{-k} < \ldots < \lambda_{-1} < 1 < \lambda_1 < \ldots < \lambda_k \tag{7.2}$$

such that $\phi^*$ acts orthogonally on $E_0$ and if $1 \leq |j| \leq k$ then
1. $\dim(E_{-j}) = \dim(E_j)$
2. $\lambda_j \lambda_{-j} = 1$
3. $\phi^*$ acts by multiplication by $\lambda_j$ on $E_j$

By construction, for all $v \in \mathrm{H}^1(S;\mathbb{R})$ and all $t \in \mathbb{R}$, $\langle v, v \rangle_{t+1} = \langle \phi^* v, \phi^* v \rangle_t$. Then given $v_0 \in E_0$ and $v_j \in E_j$, we have that for all $t \in [0,1]$ and $n \in \mathbb{Z}$,

$$\begin{aligned} \langle v_0, v_0 \rangle_{t+n} &= \langle v_0, v_0 \rangle_t, \\ \langle v_j, v_j \rangle_{t+n} &= \lambda_j^{2n} \langle v_i, v_i \rangle_t, \end{aligned} \tag{7.3}$$

Thus there is a constant $C > 0$ such that for $t \geq 0$,

$$\begin{aligned} C^{-1} \|v_0\|_0 \leq \|v_0\|_t &\leq C \|v_0\|_0, \\ C^{-1} e^{t \log(\lambda_j)} \|v_j\|_0 \leq \|v_j\|_t &\leq C e^{t \log(\lambda_j)} \|v_j\|_0. \end{aligned} \tag{7.4}$$

From Corollary 2, $0 \notin \sigma\left(\delta d \text{ on } \Lambda^1(M)/\mathrm{Ker}(d)\right)$ if and only if $E_0 = 0$.
**End of Example 1**

Example 1 shows the nicest possible behavior for $\|\cdot\|_t$. We expect that in some sense, an end of a manifold $N$ as in Corollary 3 will generally have a similar Lyapunov-type decomposition for the cohomology group $\mathrm{H}^1(S;\mathbb{R})$. We discuss the evidence for this at the end of the section. For now, we just give some consequences of having such a decomposition.

First, we give a sufficient condition for zero to not be in $\sigma\left(\delta d \text{ on } \Lambda^1(N)/\mathrm{Ker}(d)\right)$.

**Lemma 5.** *Let $V$ be a finite-dimensional real vector space with a smooth family of inner products $\{\langle\cdot,\cdot\rangle_t\}_{t\in[0,\infty)}$. Let $L^2([0,\infty);V)$ be the space of measurable maps $f:[0,\infty)\to V$ such that*

$$\|f\|^2 = \int_0^\infty \langle f(t), f(t)\rangle_t \, dt < \infty. \tag{7.5}$$

*Suppose that there are constants $a, c > 0$ such that if $s_1 \geq s_2 \geq 0$ and $v \in V$ then*

$$\|v\|_{s_1} \geq c\, e^{a(s_1-s_2)} \|v\|_{s_2}. \tag{7.6}$$

*Let $\mathcal{O}$ be the operator*

$$(\mathcal{O}f)(t) = \int_t^\infty f(s)ds. \tag{7.7}$$

*Then $\mathcal{O}$ is a bounded operator on $L^2([0,\infty);V)$.*

*Proof.* If $f \in C_0^\infty([0,\infty);V)$ then the $L^2$-norm of $\mathcal{O}f$ is given by

$$\begin{aligned} \|\mathcal{O}f\|^2 &= \int_0^\infty \langle \int_t^\infty f(s_1)ds_1, \int_t^\infty f(s_2)ds_2\rangle_t \, dt \\ &= \int_0^\infty \int_0^\infty \int_0^{\min(s_1,s_2)} \langle f(s_1), f(s_2)\rangle_t \, dt ds_1 ds_2 \\ &\leq \int_0^\infty \int_0^\infty \int_0^{\min(s_1,s_2)} \|f(s_1)\|_t \cdot \|f(s_2)\|_t \, dt ds_1 ds_2. \end{aligned} \tag{7.8}$$

Suppose that $s_1 \geq s_2 \geq s_3 \geq 0$. Then from (7.6),

(7.9)
$$\begin{aligned} \|f(s_1)\|_{s_3} \cdot \|f(s_2)\|_{s_3} &\leq c^{-1} e^{-a(s_1-s_3)} \|f(s_1)\|_{s_1} \cdot c^{-1} e^{-a(s_2-s_3)} \|f(s_2)\|_{s_2} \\ &= c^{-2}\, e^{-a(s_1-s_2)}\, e^{-2a(s_2-s_3)}\, \|f(s_1)\|_{s_1} \cdot \|f(s_2)\|_{s_2}. \end{aligned}$$

Thus if $s_1 \geq s_2$ then

(7.10)
$$\begin{aligned}\int_0^{\min(s_1,s_2)} \|f(s_1)\|_t \cdot \|f(s_2)\|_t \, dt &= \int_0^{s_2} \|f(s_1)\|_t \cdot \|f(s_2)\|_t \, dt \\ &\leq \int_0^{s_2} c^{-2} e^{-a(s_1-s_2)} e^{-2a(s_2-t)} \|f(s_1)\|_{s_1} \cdot \|f(s_2)\|_{s_2} \, dt \\ &\leq \frac{1}{2ac^2} \, e^{-a(s_1-s_2)} \, \|f(s_1)\|_{s_1} \cdot \|f(s_2)\|_{s_2}.\end{aligned}$$

In any case,

$$\int_0^{\min(s_1,s_2)} \|f(s_1)\|_t \cdot \|f(s_2)\|_t \, dt \leq \frac{1}{2ac^2} \, e^{-a|s_1-s_2|} \, \|f(s_1)\|_{s_1} \cdot \|f(s_2)\|_{s_2}. \tag{7.11}$$

Then

$$\|\mathcal{O}f\|^2 \leq \int_0^\infty \int_0^\infty \frac{1}{2ac^2} \, e^{-a|s_1-s_2|} \, \|f(s_1)\|_{s_1} \cdot \|f(s_2)\|_{s_2} \, ds_1 ds_2. \tag{7.12}$$

For $s \geq 0$, put $g(s) = \|f(s)\|_s$. Extend $g$ by zero to become an $L^2$-function on $\mathbb{R}$. Then

(7.13)
$$\int_{-\infty}^\infty \int_{-\infty}^\infty \frac{e^{-a|s_1-s_2|}}{2a} g(s_1) \, g(s_2) \, ds_1 ds_2 = \langle g, (-\partial_s^2 + a^2)^{-1} g \rangle \leq a^{-2} \int_{-\infty}^\infty g^2(s) ds.$$

The proposition follows. □

**Lemma 6.** *Let $V$ be a finite-dimensional real vector space with a smooth family of inner products $\{\langle \cdot, \cdot \rangle_t\}_{t \in [0,\infty)}$. Let $L^2([0,\infty); V)$ be the space of measurable maps $f : [0,\infty) \to V$ such that*

$$\|f\|^2 = \int_0^\infty \langle f(t), f(t) \rangle_t \, dt < \infty. \tag{7.14}$$

*Suppose that there are constants $a, c > 0$ such that if $s_1 \geq s_2 \geq 0$ and $v \in V$ then*

$$\|v\|_{s_1} \leq c \, e^{-a(s_1-s_2)} \|v\|_{s_2}. \tag{7.15}$$

*Let $\mathcal{O}'$ be the operator*

$$(\mathcal{O}'f)(t) = \int_0^t f(s) ds. \tag{7.16}$$

*Then $\mathcal{O}'$ is a bounded operator on $L^2([0,\infty); V)$.*

*Proof.* The proof is similar to that of Lemma 5. We omit the details. □

**Proposition 22.** *Let $\mathcal{U}$ contain an end of $\mathcal{M}$ as in Section 6. Let $\gamma : [0, \infty) \to \mathcal{T}_S$ be the corresponding Teichmüller ray. Let $\langle \cdot, \cdot \rangle_t$ be the inner product on $\mathrm{H}^1(S; \mathbb{R})$ coming from $\gamma(t)$. Suppose that there is a decomposition $\mathrm{H}^1(S; \mathbb{R}) = E_+ \oplus E_-$ and constants $a, c_+, c_- > 0$ such that for all $v_+ \in E_+$, $v_- \in E_-$ and $s_1 \geq s_2 \geq 0$,*

$$\|v_+\|_{s_1} \geq c_+ \, e^{a(s_1 - s_2)} \|v_+\|_{s_2} \tag{7.17}$$

*and*

$$\|v_-\|_{s_1} \leq c_- \, e^{-a(s_1 - s_2)} \|v_-\|_{s_2}. \tag{7.18}$$

*Then $\mathrm{H}^2_{(2)}(\overline{\mathcal{U}}) = 0$.*

*Proof.* From Proposition 19, we must show that $\partial_t : \Gamma(\mathrm{H}^1) \to \Gamma'(\mathrm{H}^1)$ is onto. Given $v \in \Gamma'(\mathrm{H}^1)$, write $v(t) = v_+(t) + v_-(t)$ with $v_+(t) \in E_+$ and $v_-(t) \in E_-$. Put

$$w(t) = \int_0^t v_+(s) ds - \int_t^\infty v_-(s) ds. \tag{7.19}$$

Clearly $\partial_t w = v$. By Lemmas 5 and 6, $w \in \Gamma(\mathrm{H}^1)$. □

**Corollary 4.** *Let $N$ be as in Corollary 3. Suppose that each end of $N$ is geometrically infinite and there is a decomposition of the corresponding $\mathrm{H}^1(S; \mathbb{R})$ as in the statement of Proposition 22. Then $0 \notin \sigma \left( \delta d \text{ on } \Lambda^1(N)/\mathrm{Ker}(d) \right)$.*

*Proof.* This follows from Corollary 3.2 and Proposition 22. □

We now give a sufficient condition for zero to be in $\sigma \left( \delta d \text{ on } \Lambda^1(N)/\mathrm{Ker}(d) \right)$.

**Lemma 7.** *Let $h$ be a positive smooth function on $[0, \infty)$. Suppose that there is a constant $C > 0$ such that for all $t \geq 0$,*

$$\frac{1}{C(1+t)} \leq h(t) \leq C(1+t). \tag{7.20}$$

*Put $\Gamma' = L^2(h(t)dt)$ and*

$$\Gamma = \{ f \in \Gamma' : f \text{ is absolutely continuous and } \partial_t f \in \Gamma' \}.$$

*Then $\partial_t : \Gamma \to \Gamma'$ is not onto.*

*Proof.* Put

$$g(t) = (1+t)^{-\frac{1}{2}} (\log(1+t))^{-\frac{3}{4}} h^{-\frac{1}{2}}(t). \tag{7.21}$$

Then $g \in \Gamma'$. However,

$$\int_0^t g(s) ds \geq C^{-\frac{1}{2}} \int_0^t (1+s)^{-1} (\log(1+s))^{-\frac{3}{4}} ds = 4 C^{-\frac{1}{2}} (\log(1+t))^{\frac{1}{4}}. \tag{7.22}$$

For any $T \geq 0$,

$$\int_T^\infty (\log(1+t))^{\frac{1}{2}} h(t) dt \geq \frac{1}{C} \int_T^\infty (\log(1+t))^{\frac{1}{2}} \frac{dt}{1+t} = \infty. \tag{7.23}$$

It follows that for all $c \in \mathbb{R}$, $c + \int_0^t g(s) ds$ does not lie in $L^2(h(t)dt)$ and so $g$ cannot be in the image of $\partial_t : \Gamma \to \Gamma'$. $\square$

**Proposition 23.** *Let $\mathcal{U}$ contain an end of $\mathcal{M}$ as in Section 6. Let $\gamma : [0,\infty) \to \mathcal{T}_S$ be the corresponding Teichmüller ray. Let $\langle \cdot, \cdot \rangle_t$ be the inner product on $\mathrm{H}^1(S;\mathbb{R})$ coming from $\gamma(t)$. Suppose that there is a $v \in \mathrm{H}^1(S;\mathbb{R})$ and a $C > 0$ such that for all $t \geq 0$,*

$$\frac{1}{C\sqrt{1+t}} \leq \|v\|_t \leq C\sqrt{1+t}. \tag{7.24}$$

*Then $\mathrm{H}^2_{(2)}(\overline{\mathcal{U}}) \neq 0$.*

*Proof.* By Proposition 19, we must show that $\partial_t : \Gamma(\mathrm{H}^1) \to \Gamma'(\mathrm{H}^1)$ is not onto. Putting $h(t) = \|v(t)\|^2$, this follows from Lemma 7. $\square$

**Corollary 5.** *Let $N$ be as in Corollary 3. Suppose that some end of $N$ is geometrically finite or else there is an element $v$ of the corresponding $\mathrm{H}^1(S;\mathbb{R})$ satisfying (7.24). Then $0 \in \sigma\,(\delta d \text{ on } \Lambda^1(N)/\mathrm{Ker}(d))$.*

*Proof.* This follows from Corollary 3.2 and Proposition 23. $\square$

*Remark* : Using the results of Example 1, Corollary 2 is a special case of Corollaries 4 and 5. Other examples in which the hypotheses of Corollaries 4 and 5 are satisfied are given by hyperbolic 3-manifolds with geometrically infinite ends having the same ending laminations as periodic ends.

The question arises as to how often the assumptions of Corollaries 4 and 5 hold. The qualitative behavior of the norms $\|\cdot\|_t$, as a function of $t$, is determined by the dynamics of the projected Teichmüller geodesic $\pi \circ \gamma$ on $\mathrm{Mod}_S$. Example 1 comes from the case of a closed loop on $\mathrm{Mod}_S$. Recall that as $M$ has positive injectivity radius, $\pi \circ \gamma$ lies within a compact region of $\mathrm{Mod}_S$. It seems that the dynamics of geodesics on $\mathrm{Mod}_S$ is similar to that of Riemannian geodesics on finite volume hyperbolic manifolds with cusps, in that exceptional geodesics can be constructed which have almost any desired behavior. However, one may ask if most geodesics have some uniform behavior.

The recent work of Anton Zorich is relevant here [24, 25]. Let $S$ be a closed oriented surface of genus $g \geq 2$. Instead of talking about measured geodesic laminations on $S$, we will use the equivalent language of singular foliations $\mathcal{F}$ of $S$ with an invariant transverse measure $\mu$. Zorich considers the subset $\mathcal{O}MF$ of orientable measured foliations, or equivalently, the measured foliations arising from a closed 1-form on $S$. For generic $\mathcal{F}$, the measure $\mu$ will be a unique ergodic invariant transverse measure on $\mathcal{F}$ up to scaling. Given generic $(\mathcal{F}, \mu) \in \mathcal{O}MF$, using Oseledec's theorem, Zorich constructs a certain filtration

$$0 \subset F_{-k} \subset \ldots \subset F_{-1} \subsetneq F_0 \subset F_1 \subset \ldots \subset F_k = \mathrm{H}^1(S; \mathbb{R}) \tag{7.25}$$

and positive numbers

$$\lambda_{-k} < \ldots < \lambda_{-1} < 1 < \lambda_1 < \ldots < \lambda_k \tag{7.26}$$

with $\lambda_j \lambda_{-j} = 1$, having the following property : Pick a generic point $p \in S$. Let $l$ be a half-leaf through $p$. Take a small transverse interval $I$ at $p$. Let $\{l_n\}_{n\in\mathbb{N}}$ be the segments of $l$ from $p$ to $I$, in increasing order. That is, the first return of $l$ to $I$ gives $l_1$, the second gives $l_2$, etc. For each $n \in \mathbb{N}$, close the segment $l_n$ by a short arc along $I$ joining the endpoints of $l_n$. This gives a closed loop which represents some $h_n \in \mathrm{H}_1(S; \mathbb{R})$. Pick an arbitrary Euclidean metric $\| \cdot \|$ on $\mathrm{H}_1(S; \mathbb{R})$. Then if $i > 0$ and $f_i \in F_i \backslash F_{i-1}$,

$$\limsup_{n\to\infty} \frac{\log |f_i(h_n)|}{\log \|h_n\|} = \frac{\log(\lambda_i)}{\log(\lambda_k)}. \tag{7.27}$$

Also, if $f_0 \in F_0 \backslash F_{-1}$ then

$$\limsup_{n\to\infty} \frac{\log |f_0(h_n)|}{\log \|h_n\|} = 0. \tag{7.28}$$

**Example 2 :** Consider a pseudo-Anosov diffeomorphism as in Example 1. Let $(\mathcal{F}, \mu)$ be the corresponding stable measured foliation. Note that $(\mathcal{F}, \mu)$ may not be oriented or generic. Regardless, one can see that there is a filtration (7.25) satisfying (7.27) and (7.28). In fact, it is equivalent to the decomposition (7.1), in that $F_i = F_{i-1} \oplus E_i$.
**End of Example 2**

Zorich's results are not directly applicable to our problem as we are interested in the Teichmüller rays $\gamma$ such that $\pi \circ \gamma$ is precompact, but these are not generic. Nevertheless, one can speculate on an algorithm which in "most" cases would input the end invariants of $N$ and output whether or not zero lies in the spectrum of $\sigma\left(\delta d \text{ on } \Lambda^1(N)/\mathrm{Ker}(d)\right)$. Namely, let $N$ be as in Corollary 3 and assume that all of the ends of $N$ are geometrically infinite. For each end, describe the end invariant as a measured foliation $(\mathcal{F}, \mu)$. Apply the above procedure of following a generic leaf of $\mathcal{F}$ to obtain an increasing sequence

$$F_0 \subset F_1 \subset \ldots \subset F_k = \mathrm{H}^1(S; \mathbb{R}) \tag{7.29}$$

and numbers $1 < \lambda_1 < \ldots < \lambda_k$ satisfying (7.27) and (7.28). Then zero should not be in the spectrum of $\sigma\,(\delta d \,\text{on}\, \Lambda^1(N)/\text{Ker}(d))$ if and only if for each end of $N$, $\dim(F_0) = \text{genus}(S)$.

## 8. Reduced $L^2$-cohomology

**Definition 5.** *Let $M$ be as in the beginning of Section 2. Define the relative reduced $L^2$-cohomology groups of $M$ by*

$$\overline{\mathrm{H}}^p_{(2)}(M, \partial M) = \{\omega \in \Omega^p(M) : d\omega = \delta\omega = b^*(\omega) = 0\}. \tag{8.1}$$

There is a nondegenerate pairing

$$\int_M : \overline{\mathrm{H}}^p_{(2)}(M, \partial M) \times \overline{\mathrm{H}}^{dim(M)-p}_{(2)}(M) \longrightarrow \mathbb{R}. \tag{8.2}$$

**Proposition 24.** *Let $\mathcal{U} = (0, \infty) \times S$ contain a geometrically infinite end of the model manifold $\mathcal{M}$. Suppose that the corresponding operator $\partial_t : \Gamma(\mathrm{H}^1) \to \Gamma'(\mathrm{H}^1)$ has closed image. Then $\overline{\mathrm{H}}^1_{(2)}(\mathcal{U})$ is isomorphic to a Lagrangian subspace of $\mathrm{H}^1(S; \mathbb{R})$.*

*Proof.* From Proposition 21, $\overline{\mathrm{H}}^1_{(2)}(\overline{\mathcal{U}})$ is isomorphic to a subspace of $\mathrm{H}^1(S; \mathbb{R})$. It remains to show that this subspace is Lagrangian. The pair $(\overline{\mathcal{U}}, S)$ gives a cohomology sequence

$$\ldots \longrightarrow \overline{\mathrm{H}}^1_{(2)}(\overline{\mathcal{U}}) \xrightarrow{\alpha} \mathrm{H}^1(S; \mathbb{R}) \xrightarrow{\beta} \overline{\mathrm{H}}^2_{(2)}(\overline{\mathcal{U}}, S) \longrightarrow \ldots \tag{8.3}$$

In general, this sequence will not be weakly exact without some Fredholmness assumptions. In our case, from Proposition 19, the assumption that $\partial_t$ has closed image implies that $d_1 : \Omega^1(\overline{\mathcal{U}}) \to \Omega^2(\overline{\mathcal{U}})$ is Fredholm in the sense of [12, Definition 2.1]. Then [12, Theorem 2.2] implies that (8.3) is weakly exact at $\mathrm{H}^1(S; \mathbb{R})$. As the vector spaces involved are finite-dimensional, this is the same as exactness.

Given $x \in \overline{\mathrm{H}}^1_{(2)}(\overline{\mathcal{U}})$ and $y \in \mathrm{H}^1(S; \mathbb{R})$, one can check that

$$\int_S y \cup \alpha(x) = \int_{\overline{\mathcal{U}}} \beta(y) \cup x. \tag{8.4}$$

It follows that the intersection form on $\mathrm{H}^1(S; \mathbb{R})$ vanishes when restricted to $\text{Im}(\alpha)$. Furthermore, if $y$ is perpendicular to $\text{Im}(\alpha)$ with respect to the intersection form then $y \in \text{Ker}(\beta) = \text{Im}(\alpha)$. The proposition follows. □

**Proposition 25.** *Let $N$ and $K$ be as in Corollary 3. Assume that zero does not lie in $\sigma\,(\delta d \,\text{on}\, \Lambda^1(N)/\text{Ker}(d))$. Let $L_1 \subset \mathrm{H}^1(\partial K; \mathbb{R})$ be the Lagrangian subspace*

$$\text{Im}\left(\mathrm{H}^1(K; \mathbb{R}) \longrightarrow \mathrm{H}^1(\partial K; \mathbb{R})\right).$$

*Let $L_2$ be the Lagrangian subspace of* $\mathrm{H}^1(\partial K;\mathbb{R})$ *coming from the ends of* $N$*, as in Proposition* 24*. Then there is a short exact sequence*

$$(8.5)\quad 0 \longrightarrow \mathrm{Im}\left(\mathrm{H}^1(K,\partial K;\mathbb{R}) \to \mathrm{H}^1(K;\mathbb{R})\right) \longrightarrow \overline{\mathrm{H}}^1_{(2)}(N) \longrightarrow L_1 \cap L_2 \to 0.$$

*Proof.* By Corollary 3.2, each end of $N$ is geometrically infinite and the corresponding operator $\partial_t : \Gamma(\mathrm{H}^1) \to \Gamma'(\mathrm{H}^1)$ has closed image. Let $V$ be the closure of a union of open sets $(0,\infty) \times S_i$ containing the ends of $N$. Take $K = \overline{N - V}$. There is a Mayer-Vietoris sequence

$$(8.6)\quad \begin{aligned} \ldots \longrightarrow &\overline{\mathrm{H}}^1_{(2)}(N) \longrightarrow \mathrm{H}^1(K;\mathbb{R}) \oplus \overline{\mathrm{H}}^1_{(2)}(V) \longrightarrow \mathrm{H}^1(\partial K;\mathbb{R}) \longrightarrow \\ &\overline{\mathrm{H}}^2_{(2)}(N) \longrightarrow \mathrm{H}^2(K;\mathbb{R}) \oplus \overline{\mathrm{H}}^2_{(2)}(V) \longrightarrow \mathrm{H}^2(\partial K;\mathbb{R}) \longrightarrow \ldots \end{aligned}$$

Again, this sequence will not be weakly exact in full generality. However, in our case $d_1 : \Omega^1(V) \to \Omega^2(V)$ is Fredholm. Along with the fact that the differentials $d : \Omega^*(S) \to \Omega^{*+1}(S)$ are Fredholm, [12, Theorem 2.2] implies that (8.6) is weakly exact at the terms from $\mathrm{H}^1(K;\mathbb{R}) \oplus \overline{\mathrm{H}}^1_{(2)}(V)$ to $\mathrm{H}^2(K;\mathbb{R}) \oplus \overline{\mathrm{H}}^2_{(2)}(V)$. Again, as the vector spaces are finite-dimensional, the sequence will actually be exact at these terms. By Proposition 20, $\overline{\mathrm{H}}^2_{(2)}(V) = 0$. Dualizing (8.6) gives a sequence

(8.7)
$$\begin{aligned} \ldots \longrightarrow &\mathrm{H}^0(\partial K;\mathbb{R}) \longrightarrow \mathrm{H}^1(K,\partial K;\mathbb{R}) \longrightarrow \overline{\mathrm{H}}^1_{(2)}(N) \longrightarrow \\ &\mathrm{H}^1(\partial K;\mathbb{R}) \longrightarrow \mathrm{H}^2(K,\partial K;\mathbb{R}) \oplus \overline{\mathrm{H}}^2_{(2)}(V,\partial K) \longrightarrow \overline{\mathrm{H}}^2_{(2)}(N) \longrightarrow \ldots \end{aligned}$$

which is exact at the terms from $\mathrm{H}^1(K,\partial K;\mathbb{R})$ to $\mathrm{H}^2(K,\partial K;\mathbb{R}) \oplus \overline{\mathrm{H}}^2_{(2)}(V,\partial K)$. This gives the short exact sequence

$$(8.8)\quad \begin{aligned} 0 \longrightarrow &\mathrm{Coker}\left(\mathrm{H}^0(\partial K;\mathbb{R}) \longrightarrow \mathrm{H}^1(K,\partial K;\mathbb{R})\right) \longrightarrow \overline{\mathrm{H}}^1_{(2)}(N) \longrightarrow \\ &\mathrm{Ker}\left(\mathrm{H}^1(\partial K;\mathbb{R}) \longrightarrow \mathrm{H}^2(K,\partial K;\mathbb{R}) \oplus \overline{\mathrm{H}}^2_{(2)}(V,\partial K)\right) \longrightarrow 0. \end{aligned}$$

From the exact cohomology sequence of the pair $(K,\partial K)$,

$$(8.9)\quad \mathrm{Coker}\left(\mathrm{H}^0(\partial K;\mathbb{R}) \longrightarrow \mathrm{H}^1(K,\partial K;\mathbb{R})\right) \cong \mathrm{Im}\left(\mathrm{H}^1(K,\partial K;\mathbb{R}) \to \mathrm{H}^1(K;\mathbb{R})\right)$$

and

(8.10)
$$\mathrm{Ker}\left(\mathrm{H}^1(\partial K;\mathbb{R}) \longrightarrow \mathrm{H}^2(K,\partial K;\mathbb{R})\right) \cong \mathrm{Im}\left(\mathrm{H}^1(K;\mathbb{R}) \longrightarrow \mathrm{H}^1(\partial K;\mathbb{R})\right) = L_1.$$

Thus

$$(8.11)\quad \begin{aligned} &\mathrm{Ker}\left(\mathrm{H}^1(\partial K;\mathbb{R}) \longrightarrow \mathrm{H}^2(K,\partial K;\mathbb{R}) \oplus \overline{\mathrm{H}}^2_{(2)}(V,\partial K)\right) = \\ &L_1 \cap \mathrm{Ker}\left(\mathrm{H}^1(\partial K;\mathbb{R}) \longrightarrow \overline{\mathrm{H}}^2_{(2)}(V,\partial K)\right). \end{aligned}$$

Identifying $\overline{\mathrm{H}}^1_{(2)}(V)$ with the subspace $L_2$ of $\mathrm{H}^1(\partial K;\mathbb{R})$, the pairing (8.2) gives

$$\overline{\mathrm{H}}^2_{(2)}(V,\partial K) \cong \left(\overline{\mathrm{H}}^1_{(2)}(V)\right)^* \cong L_2^*. \tag{8.12}$$

The map $A : \mathrm{H}^1(\partial K;\mathbb{R}) \longrightarrow \overline{\mathrm{H}}^2_{(2)}(V,\partial K) \cong L_2^*$ is given explicitly by

$$(A(h))\,(l) = \int_S h \cup l \tag{8.13}$$

for all $h \in \mathrm{H}^1(\partial K;\mathbb{R})$ and $l \in L_2$. As $L_2$ is Lagrangian,

$$\mathrm{Ker}\left(\mathrm{H}^1(\partial K;\mathbb{R}) \longrightarrow \overline{\mathrm{H}}^2_{(2)}(V,\partial K)\right) = L_2. \tag{8.14}$$

The proposition now follows from equations (8.8), (8.9), (8.11) and (8.14). □

**Example 3 :** Let $M$ be as in Example 1, with $E_0 = 0$. With respect to the diffeomorphism $M = \mathbb{R} \times S$, take $K = [-1,1] \times S$. Then $M$ certainly satisfies the hypotheses of Proposition 25. We have $\partial K = S \amalg S$, with the Lagrangian subspace $L_1$ being the diagonal in $\mathrm{H}^1(K;\mathbb{R}) = \mathrm{H}^1(S;\mathbb{R}) \oplus \mathrm{H}^1(S;\mathbb{R})$. As $L_2 = \left(\bigoplus_{i=1}^k E_i\right) \bigoplus \left(\bigoplus_{i=1}^k E_{-i}\right)$, we have $L_1 \cap L_2 = 0$. Then Proposition 25 gives $\overline{\mathrm{H}}^1_{(2)}(M) = 0$. Of course, this is consistent with Proposition 16.1.

Now let $Z$ be the subset $[0,\infty) \times S$ of $M$. Perturb the metric on $Z$ to make it a product near $\{0\} \times S$. Let $N$ be the double of $Z$. Again, $N$ is diffeomorphic to $\mathbb{R} \times S$. Take $K = [-1,1] \times S$. Then $N$ also satisfies the hypotheses of Proposition 25. In this case, $L_2 = \left(\bigoplus_{i=1}^k E_{-i}\right) \bigoplus \left(\bigoplus_{i=1}^k E_{-i}\right)$. Thus $L_1 \cap L_2 = L_2$. Proposition 25 gives $\dim\left(\overline{\mathrm{H}}^1_{(2)}(N)\right) = g$, the genus of $S$. This shows that in the setting of Proposition 25, $\overline{\mathrm{H}}^1_{(2)}(N)$ depends on the end invariants of $N$ and not just on the topological type of $K$.

**End of Example 3**

Department of Mathematics, University of Michigan, Ann Arbor, MI 48109, USA
*E-mail address*: lott@math.lsa.umich.edu